\documentclass{article}
\title{Stabilizing Parabolic Vector Bundles}
\author{Andy Ramirez-Cote}
\date{November 2026}
\usepackage[english]{babel}

\usepackage[letterpaper,top=2cm,bottom=2cm,left=3cm,right=3cm,marginparwidth=1.75cm]{geometry}

\usepackage{amsmath}
\usepackage{amsfonts}
\usepackage{amsthm}
\usepackage{amssymb}
\usepackage{mathrsfs}
\usepackage[T1]{fontenc}%
\usepackage{microtype}
\usepackage{calligra}
\usepackage{graphicx}
\usepackage{tikz-cd}
\usepackage{csquotes}
\usepackage{xparse}
\usepackage[style=alphabetic, backend=biber]{biblatex}
\usepackage[colorlinks=true, allcolors=blue]{hyperref}
\DeclareDocumentCommand{\newmathcommand}{mO{0}m}{%
  \expandafter\NewCommandCopy\csname old\string#1\endcsname{#1}%
  \expandafter\newcommand\csname new\string#1\endcsname[#2]{#3}
  \DeclareRobustCommand#1{%
    \ifmmode
      \expandafter\let\expandafter\next\csname new\string#1\endcsname
    \else
      \expandafter\let\expandafter\next\csname old\string#1\endcsname
    \fi
    \next
  }%
}

\theoremstyle{definition}
\newtheorem{theorem}{Theorem}[section]
\newtheorem{definition}[theorem]{Definition}
\newtheorem{lemma}[theorem]{Lemma}
\newtheorem{proposition}[theorem]{Proposition}

\newtheorem{theorem-definition}[theorem]{Theorem-Definition}
\newtheorem{corollary}[theorem]{Corollary}

\newtheorem{remark}[theorem]{Remark}
\newtheorem{notation}[theorem]{Notation}

\DeclareMathOperator{\Hom}{Hom}
\DeclareMathOperator{\ihom}{\mathcal{H}\mkern-1.5mu\mathit{om}}
\DeclareMathOperator{\End}{\mathcal{E}\mkern-1.5mu\mathit{nd}}
\DeclareMathOperator{\Quot}{\mathcal{Q}\mkern-1.5mu\mathit{uot}}

\DeclareMathOperator{\im}{im}
\DeclareMathOperator{\gr}{gr}

\DeclareMathOperator{\GL}{GL}

\newcommand{\pa}[1]{{#1}_{\star}}

\newcommand{\rk}{\operatorname{rk}}
\newcommand{\pardeg}{\operatorname{par-deg}}
\newmathcommand{\len}{\operatorname{len}}
\newmathcommand{\pmu}{\pa{\mu}}
\ExplSyntaxOn
\NewDocumentCommand{\aAt}{m}
{
    \operatorname{At}\sb{C}(#1)
}
\NewDocumentCommand{\At}{m}
{
    \operatorname{At}\sb{(C,D)}(
        \clist_item:nn {#1}{1}
        \int_compare:nNnT {\clist_count:n {#1}} > {1}
            {,\clist_item:nn {#1}{2}\sb{\bullet}}
    )
}
\NewDocumentCommand{\pAt}{m}
{
    \operatorname{At}\sb{(C,D)}(
        \pa{\clist_item:nn {#1}{1}}
        \int_compare:nNnT {\clist_count:n {#1}} > {1}
            {,\clist_item:nn {#1}{2}\sb{\bullet}}
    )
}

\ExplSyntaxOff
\newmathcommand{\parvect}{\pa{\operatorname{Vect}}(C,D)}
\newmathcommand{\coparvect}{\widehat{\pa{\operatorname{Vect}}}(C,D)}

\newmathcommand{\fdef}{$Spf(\widehat{\mathcal{O}}_{X,x}/\mathcal{I_{\nabla}})$}

\begin{document}

\maketitle

\begin{abstract}
For an irreducible parabolic flat vector bundle of rank $r < 2\sqrt{g}$ on a genus $g$ curve,
we prove that an analytically general isomonodromic deformation is stable. This refines
Landesman--Litt's semistability result for bundles of rank $< 2\sqrt{g+1}$.

The proof applies the non-generically globally generated (non-GGG) lemma to
establish infinitesimal obstructions to extending destabilizing subbundles of
semistable bundles.
\end{abstract}

\section{Introduction}
We fix the following notation for the rest of this article. $(C,D)$ will always 
be an algebraic curve of genus $g \geq 2$ together with a reduced divisor 
$D \in \text{Div}(C)$, and $(V,\nabla)$ always denotes some flat vector bundle on $C$. 
This paper is about the stability properties of generic deformation of a triple $(C,V,\nabla)$ preserving 
the associated monodromy representation.
Our main result establishes that for an irreducible parabolic flat vector bundle $(V,\nabla)$ of rank $r$ on a genus $g$ curve $C$,
an analytically general isomonodromic deformation is stable, provided the rank satisfies $r < 2\sqrt{g}$.
\begin{theorem}\label{thm:intro-stability}
Let $(C,D)$ be a smooth projective curve of genus $g$ with marked points, and let
$(\pa{V},\nabla)$ be the regular flat parabolic bundle associated to an irreducible local system
of rank $r$ on $C \setminus D$. If $r < 2\sqrt{g}$, then an analytically general isomonodromic deformation
of $\pa{V}$ is stable, and each of its nonzero proper parabolic subbundles $\pa{V'}$ satisfies
\[\pardeg \pa{V'} \leq -\frac{g - \lfloor r^{2}/4 \rfloor}{r}.\]
\end{theorem}

Here, `analytically general' is meant in the sense of \cite[Definition 1.2.3]{landesman2023geometric}:
locally on $\mathcal{M}_{g,n}$, outside a nowhere dense closed analytic subset. An isomonodromic deformation
preserves the monodromy representation of the flat connection (see Section~\ref{sec:isomonodromy}
for precise definitions), and parabolic structure together with its associated degree are
reviewed in Section~\ref{sec:parabolic}.

\subsection{Previous work}
Biswas, Heu and Hurtubise asked whether a flat vector bundle with regular singularities always admits a semistable isomonodromic deformation \cite[p.~123]{BiswasHeuHurtubise2016}. Landesman and Litt showed that this can fail even when the monodromy is irreducible \cite[Corollary 1.3.3]{landesman2023geometric}, and that it nonetheless holds for irreducible monodromy when the rank is less than $2\sqrt{g+1}$ \cite[Corollary 6.1.2]{landesman2023geometric}. Theorem~\ref{thm:intro-stability} refines this semistability result to full stability.

For vector bundles of rank $2$ more is known. Heu proved that a general isomonodromic deformation is stable when $g\geq 2$ \cite{Heu2009}, and Biswas, Heu and Hurtubise proved the stronger result that it does not admit a nilpotent Higgs field \cite[Theorem 1.2]{BiswasHeuHurtubise2017}, which implies stability by a classical result of Laumon \cite[Proposition 3.5]{Laumon1988}.

\subsection{Structure of the paper}
We review the notion of parabolic bundles in Section 2 following \cite{landesman2023geometric}, and for the complete exposition we refer the reader to \cite{Yokogawa1995}.
Section 3 reviews the Atiyah bundle and its variants for filtered and parabolic bundles and their relation to deformation theory. Some explicit verifications are written down for results which are folklore but, to the best of the author's knowledge, hadn't been written down publicly, such as Proposition \ref{prop:deformation-atiyah}. The most important result for later sections is Corollary \ref{cor:filtration-obstruction} which describes the obstruction to extending a  filtration to a first-order deformation of a parabolic vector bundle together with the underlying curve.
Section 4 covers isomonodromic deformations of $(\pa{V},\nabla)$ and reduces our problem to a purely algebraic one about first-order deformations of a filtered bundle on a single curve. We start by reminding the reader of some facts concerning isomonodromic deformations in the analytic setting, where the classical construction naturally lives over the universal cover of the analytification of $\mathcal{M}_{g,n}$. We then introduce formal and first-order isomonodromic deformations, and show in Theorem~\ref{serenity} that a parabolic subbundle extends along the analytic isomonodromic deformation if and only if it extends along the formal one, in which case it extends to first order in every direction. Along the way we show that the Quot functors of the universal bundle over the de Rham moduli space are represented by proper analytic spaces (Theorem~\ref{thm:quot}).
Section 5 concerns the non-generically globally generated
(non-GGG) lemma \cite[Lemma 6.2.3, Proposition 6.3.1]{landesman2023geometric} of Landesman and Litt, which constrains bundles via numerical conditions on filtrations. We give the proof of the non-GGG lemma for the reader's convenience, and derive some corollaries more easily applicable for our purposes. Section 6 consists in its entirety of the proof of Theorem \ref{thm:intro-stability}.

\subsection{Open problems}
The question of sharp bounds remains open. Landesman and Litt construct examples
of vector bundles of rank at most exponential in $g$ all of whose isomonodromic deformations fail
to be semistable \cite[Theorem 1.3.2, Corollary 1.3.3, Example 5.1.6]{landesman2023geometric}.
The gap between the $2\sqrt{g}$ bound proven 
here and these exponential examples leaves significant room for improvement in understanding the true thresholds for semistability{} and stability. The true thresholds for stability under isomonodromic deformation are, in the author's opinion, much higher than the bound proved here would suggest, whose asymptotic form appears to be due to the use of the non-GGG lemma and not anything intrinsic to the problem. As such it remains open to find approaches which could provide improved asymptotics. Finally, it remains unclear whether analogous results hold for very stability. If so, the isomonodromy flow would make parabolic vector bundles "as stable as possible". Nothing is known for this in rank $3$ or above.

\subsection{Acknowledgements}
I am deeply grateful to my advisor Daniel Litt for his explanations, his advice, and more hours of his time than I can count. I thank Aaron Landesman for pointing out that a destabilizing subbundle of a semistable bundle can be assumed to be stable, which led to various technical simplifications.
Claude provided valuable editorial assistance, helped find
appropriate references and verify references against their source, and suggested the reformulation of the statement of Theorem \ref{thm:stability} which became Theorem \ref{thm:intro-stability} and Corollary \ref{cor:stability}.

\section{Parabolic vector bundles}\label{sec:parabolic}
We begin by recalling some relevant definitions 
and results concerning parabolic vector bundles for the reader's convenience. Our exposition
is mainly adapted from \cite[\S2]{landesman2023geometric} and \cite[\S1]{Yokogawa1995}, but for another 
perspective on the basics on parabolic bundles see the classic reference
\cite{MehtaSeshadri1980}.

Throughout this section, we use $C$ to denote a smooth projective curve over $\mathbb{C}$
and $D = x_{1} + \cdots + x_{n}$ a reduced effective divisor on $C$. We write 
$\operatorname{Vect}(C)$ for the category of vector bundles on $C$ and 
$\operatorname{Shv}(C)$ for the category of sheaves of quasi-coherent $\mathcal{O}_{C}$-modules.

\subsection{Parabolic and coparabolic bundles}

\begin{definition}
    Let $V$ be a vector bundle on $C$ and $D=x_{1}+\cdots+x_{n}$ a reduced divisor.
    For any $j$, let $V_{x_{j}}$ be the restriction of $V$ to the point $x_{j}$.
    A quasi-parabolic structure consists of a set of $n$ strictly decreasing filtrations
    \begin{equation}
        V_{x_{j}} = V^{1}_{x_{j}} \supsetneq V^{2}_{x_{j}} \supsetneq \cdots \supsetneq V^{n_{j}+1}_{x_{j}} = 0.
    \end{equation}
    For convenience we abbreviate $V^{i}_{j} := V^{i}_{x_{j}}$.
    
    A parabolic structure on $V$ is a quasi-parabolic structure together with, 
    for each $j \in \{1,\ldots,n\}$, a sequence of real numbers 
    $\{\alpha^{i}_{j}\}^{n_{j}}_{i=1}$, which we call the parabolic weights, satisfying
    \begin{equation}
        0 \leq \alpha^{1}_{j} < \alpha^{2}_{j} < \cdots < \alpha^{n_{j}}_{j} < \alpha^{n_{j}+1}_{j} = 1.
    \end{equation}
    A parabolic vector bundle $\pa{V}$ is a triple $(V, \{V^{i}_{j}\}^{n_{j}+1}_{i=1}, \{\alpha^{i}_{j}\}^{n_{j}}_{i=1})$. 
    We define $\len(\pa{V},j) := n_{j}$ to be the number of nontrivial 
    steps in the filtration at $x_{j}$.
\end{definition}

\begin{definition}\label{def:induced-parabolic}
    Let $\pa{V}$ be a parabolic bundle on $(C,D)$ and let $V' \subset V$ be a subbundle 
    of its underlying bundle. Let 
    \[V'^{i}_{j} := V'_{x_{j}}\cap V^{i}_{j} \] 
    and
    \[ \beta^{i}_{j} := \max_{1 \leq k \leq n_{j}}\{\alpha^{k}_{j} \mid V'^{i}_{j} = V'_{x_{j}} \cap V^{k}_{j}\}.\]
    The induced quasi-parabolic structure on $V'$ is given by the set $\{V'^{i}_{j}\}$. 
    The induced parabolic structure on $V'$ is defined as $\pa{V}' := (V',\{V'^{i}_{j}\},\{\beta^{i}_{j}\})$. 
\end{definition}

\begin{definition}
     Let $\pa{V}$ be a parabolic bundle on $(C,D)$ and let $\pi: V \twoheadrightarrow Q$ 
     be a quotient of its underlying bundle.
     Let 
     \[Q^{i}_{j} := \pi(V^{i}_{j}) \] 
     and
    \[ \beta^{i}_{j} := \max_{1 \leq k \leq n_{j}}\{\alpha^{k}_{j} \mid Q^{i}_{j} = \pi(V^{k}_{j})\}.\]
    The induced quasi-parabolic structure on $Q$ is given by the set $\{Q^{i}_{j}\}$. 
    The induced parabolic weight on $Q^{i}_{j}$ is $\beta^{i}_{j}$.
    
    The induced parabolic structure on $Q$ is defined as $\pa{Q} := (Q,\{Q^{i}_{j}\},\{\beta^{i}_{j}\})$.
\end{definition}

\begin{definition}
    Given two parabolic bundles $\pa{V} = (V,\{V^{i}_{j}\},\{\alpha^{i}_{j}\})$ and 
    $\pa{V}' = (V',\{V'^{i}_{j}\},\{\beta^{i}_{j}\})$ on $C$, a morphism $\varphi: V \to V'$
    of vector bundles is a morphism of parabolic vector bundles $\pa{V} \to \pa{V}'$
    if for every $1 \leq j \leq n$ and $1 \leq i \leq n_{j}$ 
    \begin{equation}\label{morphismdef} 
        \alpha^{i}_{j} \leq \max_{1\leq k \leq \len(\pa{V}',j)+1}
        \{\beta^{k}_{j} \mid \varphi_{x_{j}}(V^{i}_{j}) \subseteq V'^{k}_{j} \}.
    \end{equation}
    Parabolic morphisms from $\pa{V}$ to $\pa{V}'$ form a subspace of $\Hom(V,V')$ which we'll denote $\pa{\Hom}(V,V')$.
    For parabolic morphisms $\varphi: \pa{V} \to \pa{V}'$
    and $\psi: \pa{V}' \to \pa{V}''$, with weights $\gamma^{\ell}_{j}$ on $\pa{V}''$,
    \[
        \begin{aligned}
            \alpha^{i}_{j}
            &\leq \max_{\varphi_{x_{j}}(V^{i}_{j}) \subseteq V'^{k}_{j}} \beta^{k}_{j}
            \leq \max_{\varphi_{x_{j}}(V^{i}_{j}) \subseteq V'^{k}_{j}}
                 \max_{\psi_{x_{j}}(V'^{k}_{j}) \subseteq V''^{\ell}_{j}} \gamma^{\ell}_{j}\\
            &\leq \max_{(\psi\circ\varphi)_{x_{j}}(V^{i}_{j}) \subseteq V''^{\ell}_{j}}
                 \gamma^{\ell}_{j}.
        \end{aligned}
    \]
    We define the category $\parvect$ to be the category whose objects 
    are parabolic vector bundles on $C$ whose quasi-parabolic structure is supported on $D$ and whose morphisms are as just defined.
\end{definition}

\begin{definition}\label{def:coparabolic}
    Let $V$ be a vector bundle on $C$ and $(\{V^{k}_{j}\}^{m_{j}}_{k=0})$ a quasi-parabolic structure on $V$. 
    A coparabolic structure on $V$ consists of a quasi-parabolic structure together with, for each $j \in \{1, \cdots, n\}$, a sequence of real numbers $\{\beta^{i}_{j}\}^{m_{j}}_{i=1}$ satisfying
    \[0 < \beta^{1}_{j} < \cdots < \beta^{m_{j}}_{j} \leq 1\]
    which we refer to as the coparabolic weights of the coparabolic structure. 
    A coparabolic vector bundle is a triple $\widehat{\pa{V}} = (V, \{V^{k}_{j}\}^{m_{j}}_{k=0}, \{\beta^{k}_{j}\}^{m_{j}}_{k=1})$ consisting of a vector bundle together with a coparabolic structure.
    A morphism $\varphi: V \to V'$ of vector bundles is a morphism of coparabolic 
    vector bundles $\widehat{\pa{V}} \to \widehat{\pa{V}'}$ if for every $1 \leq j \leq n$ 
    and $0 \leq k < m_{j}$
    \[
        \beta^{k+1}_{j} \leq \min_{1 \leq \ell \leq \len(\widehat{\pa{V}'},j)} 
        \{\gamma^{\ell}_{j} \mid \varphi_{x_{j}}(V^{k}_{j}) \not\subseteq V'^{\ell}_{j}\}
    \] 
    Just like before a simple verification shows that the identity is always coparabolic and compositions of coparabolic morphisms are coparabolic.
    We thus analogously define the category $\coparvect$ to have as objects coparabolic vector bundles on $C$ whose quasi-parabolic structure is supported on $D$ and coparabolic morphisms as morphisms.
\end{definition}

\begin{definition}\label{def:corresponding-coparabolic}
    Let $\pa{V} = (V, \{V^{i}_{j}\}, \{\alpha^{i}_{j}\}) \in \operatorname{Ob}(\parvect)$. 
    Set $V' := \ker(V \to \bigoplus_{\alpha^{1}_{j}=0} V_{x_{j}}/V^{2}_{j})$
    and let $\rho_{j}: V'_{x_{j}} \to V_{x_{j}}$ be the induced map on fibers.
    The corresponding coparabolic bundle is
    $\widehat{\pa{V}} := (V',\{V'^{k}_{j}\}^{n_{j}}_{k=0},\{\beta^{k}_{j}\}^{n_{j}}_{k=1})$,
    where $V'^{n_{j}}_{j}:=0$ and, for $1 \leq k \leq n_{j}$,
    \[
        (V'^{k-1}_{j},\beta^{k}_{j}) :=
        \begin{cases}
            (\rho_{j}^{-1}(V^{k}_{j}),\alpha^{k}_{j}),
                & \alpha^{1}_{j}>0,\\
            (\rho_{j}^{-1}(V^{k+1}_{j}),\alpha^{k+1}_{j}),
                & \alpha^{1}_{j}=0.
        \end{cases}
    \]
\end{definition}

\begin{proposition}\label{prop:parabolic-coparabolic-equivalence}
    Sending $\pa{V}$ to $\widehat{\pa{V}}$, and a parabolic morphism $\varphi: \pa{V} \to \pa{W}$ to its restriction
    $\varphi|_{V'}$, defines an equivalence of categories $\parvect \to \coparvect$.
\end{proposition}

\begin{proof}
    Definition~\ref{def:corresponding-coparabolic} is the elementary modification of $V$ along its steps of
    positive weight, which parabolic morphisms preserve, with the flags pulled back along $\rho_{j}$ and
    completed by $\ker \rho_{j}$, and with the weights permuted cyclically, so that $\alpha^{1}_{j} = 0$ becomes
    $\beta^{n_{j}}_{j} = 1$. Elementary modifications are functorial and the square formed by $\rho_{j}$, its
    analogue for $\pa{W}$ and the fibers of $\varphi$ and $\varphi|_{V'}$ commutes, so the flags are carried
    along and, once the indices are permuted, the condition~\eqref{morphismdef} on $\varphi$ is the
    coparabolic condition on $\varphi|_{V'}$. Each step is reversible, the inverse functor being the
    elementary modification of $W(D)$ along $W^{m_{j}-1}_{j}(D)$ wherever $\beta^{m_{j}}_{j} = 1$, with the image
    flags and the inverse permutation.
\end{proof}

\begin{definition}\label{def:parabolic-internal-hom}
    Let $\pa{V},\pa{V}' \in \operatorname{Ob}(\parvect)$. We define the parabolic internal hom
    $\pa{\ihom}(V,V')$ to be the subpresheaf of $\ihom(V,V')$ whose sections over 
    an open set $U \subseteq C$ consist of morphisms in $\pa{\Hom}(V|_{U},V'|_{U})$. We abbreviate $\pa{\ihom}(V,V)$ by $\pa{\End}(V)$.
\end{definition}

\begin{proposition}\label{prop:parabolic-internal-hom-bundle}
    Let $\pa{V}$ and $\pa{V}' \in \operatorname{Ob}(\parvect)$ be parabolic vector bundles on $C$. 
    Then $\pa{\ihom}(V,V')$ is a vector bundle on $C$ of rank $\rk V \cdot \rk V'$.
\end{proposition}

\begin{proof}
    We first verify that $\pa{\ihom}(V,V')$ is a sheaf. Let $U \subseteq C$ be open with open cover $\{U_{i}\}_{i \in I}$, 
    and let $\{\varphi_{i}\}$ be a compatible family of sections $\varphi_{i} \in \pa{\ihom}(V,V')(U_{i})$. 
    Since $\pa{\ihom}(V,V')$ is a subpresheaf of $\ihom(V,V')$, there exists a unique 
    $\varphi \in \ihom(V,V')(U)$ restricting to each $\varphi_{i}$. 
    
    For any $x \in D \cap U$, choosing any $i$ with $x \in U_{i}$ we find that $\varphi_{x} = (\varphi|_{U_{i}})_{x} = (\varphi_{i})_{x}$ 
    satisfies the parabolic morphism condition at $x$ since $\varphi_{i}$ does and thus $\varphi \in \pa{\ihom}(V,V')(U)$.

    As to local freeness, note that away from $D$, we have $\pa{\ihom}(V,V') = \ihom(V,V')$, which is locally free, and so at each $x_{j} \in D$ we choose a sufficiently small neighborhood 
    $U_{j}$ such that $V|_{U_{j}}$ and $V'|_{U_{j}}$ are trivial and $U_{j} \cap D = \{x_{j}\}$. 
    The parabolic condition at $x_{j}$ defines a linear subspace $P_{j}$ of
    $\Hom(V_{x_{j}}, V'_{x_{j}})$, thus $\pa{\ihom}(V,V')|_{U_{j}}$ is the kernel of $\ihom(V,V')|_{U_{j}} \to \Hom(V_{x_{j}}, V'_{x_{j}})/P_{j}$, hence locally free.
\end{proof}

We note in passing that every ordinary vector bundle can be thought of as a trivially parabolic vector bundle or trivially coparabolic bundle in a sense which we now make precise:
\begin{definition}\label{def:trivial-parabolic-structure}
    Let $V$ be a vector bundle. The trivial parabolic (resp. coparabolic) structure on $V$ consists of the quasi-parabolic structure given by the filtration $V_{x_{j}} = V^{1}_{j} \supset V^{2}_{j} = 0$ at each $x_j$ together with its unique weight $\alpha^{1}_{j}=0$ (resp. $\beta^{1}_{j}=1$). 
    
    Endowing a vector bundle with its trivial parabolic (resp. coparabolic) structure gives a functor 
    $(-)_{\mathrm{triv}}:\operatorname{Vect}(C) \to \parvect$ (resp. $\widehat{(-)_{\mathrm{triv}}}:\operatorname{Vect}(C) \to \coparvect$) left adjoint to the
    forgetful functor $\parvect \to \operatorname{Vect}(C)$ (resp. $\coparvect \to \operatorname{Vect}(C)$).
\end{definition}

\begin{definition}\label{def:parabolic-dual}
    Let $\pa{V} = (V,\{V^{i}_{j}\},\{\alpha^{i}_{j}\})$ be a parabolic vector bundle. The parabolic dual of $\pa{V}$ is
    the coparabolic bundle
    \[
        \widehat{\pa{V}^{\vee}} := \big(V^{\vee}(-D),\{\operatorname{ann}(V^{n_{j}+1-k}_{j})\}^{n_{j}}_{k=0},\{1-\alpha^{n_{j}+1-k}_{j}\}^{n_{j}}_{k=1}\big),
    \]
    with annihilators taken in $V^{\vee}(-D)_{x_{j}}$, and we call its corresponding parabolic bundle $\pa{V}^{\vee}$
    (Proposition~\ref{prop:parabolic-coparabolic-equivalence}) the naive dual of $\pa{V}$. Since the parabolic dual is in most circumstances more natural to work with, we will simply refer to it as the dual of a parabolic bundle unless this could cause confusion.
\end{definition}

\begin{proposition}\label{prop:parabolic-coparabolic-duality}
    Parabolic dualization is an antiequivalence $\parvect \to \coparvect$.
\end{proposition}

\begin{proof}
    The construction is built from invertible and functorial operations: Starting with the underlying bundles,
    $V \mapsto V^{\vee}(-D)$ is an antiequivalence of $\operatorname{Vect}(C)$ which is its own quasi-inverse. Next, taking annihilators reverses
    inclusion since $\operatorname{ann}(U)$ is the dual of $V_{x_{j}}/U$ and thus restricting to fibers we see
    that a morphism $f \colon V \to W$ maps $(V^{i}_{j},\alpha^{i}_{j})$ into $(W^{k}_{j},\beta^{k}_{j})$ if and only if $f^{\vee}$ maps
    $(\operatorname{ann}(W^{k}_{j}),1-\beta^{k}_{j})$ into $(\operatorname{ann}(V^{i}_{j}),1-\alpha^{i}_{j})$. If $f^{\vee}$ violated the coparabolic condition of
    Definition~\ref{def:coparabolic}, some $\operatorname{ann}(W^{k}_{j})$ would not be mapped into
    $\operatorname{ann}(V^{i}_{j})$ even though $\operatorname{ann}(V^{i}_{j})$ has smaller weight than
    $\operatorname{ann}(W^{k-1}_{j})$, and then $f$ would not map $V^{i}_{j}$ into $W^{k}_{j}$ even though
    $W^{k-1}_{j}$ has smaller weight than $V^{i}_{j}$, in violation of~\eqref{morphismdef}. The converse
    holds by the same logic in reverse and so dualization is fully faithful, with the same formulas applied to a
    coparabolic bundle giving a quasi-inverse.
\end{proof}

\subsection{Parabolic sheaves}

\begin{definition}
    Let $\mathbb{R}^{\mathrm{op}}$ be the opposite of the category associated to $\mathbb{R}$ 
    with its usual ordering. Explicitly, $\operatorname{Ob}(\mathbb{R}^{\mathrm{op}}) = \mathbb{R}$ and 
    $\Hom(\alpha,\beta) = \left\{\begin{smallmatrix}
        \emptyset &\text{ if } \alpha < \beta \\
        \{*\} &\text{ if } \alpha \geq \beta
    \end{smallmatrix}\right.$.
    
    A parabolic sheaf $\pa{\mathcal{F}}$ on $(C,D)$ is a functor
    $\pa{\mathcal{F}}:\mathbb{R}^{\mathrm{op}} \to \operatorname{Shv}(C)$ 
    together with a natural isomorphism $\pa{\mathcal{F}}[(-)+1] \simeq \pa{\mathcal{F}}[(-)]\otimes \mathcal{O}_{C}(-D)$ such that,
    for every $\alpha \in \mathbb{R}$, the square
    \[
    \begin{tikzcd}[column sep=large]
        {\pa{\mathcal{F}}[\alpha+1]} \arrow[r, "\sim"] \arrow[d]
            & {\pa{\mathcal{F}}[\alpha] \otimes \mathcal{O}_{C}(-D)} \arrow[d] \\
        {\pa{\mathcal{F}}[\alpha]} \arrow[r, "\sim"]
            & {\pa{\mathcal{F}}[\alpha] \otimes \mathcal{O}_{C}}
    \end{tikzcd}
    \]
    commutes. We denote the category of parabolic sheaves on $(C,D)$ by $\pa{\operatorname{Shv}}(C,D)$.
\end{definition}
\begin{remark}
Given any parabolic sheaf, restriction to the half-open interval
$[0,1)$ (resp. $(0,1]$) gives a functor $[0,1)^{\mathrm{op}} \to \operatorname{Shv}(C)$
(resp. $(0,1]^{\mathrm{op}} \to \operatorname{Shv}(C)$). Conversely, any left-continuous (resp. right-continuous) functor
$\pa{\mathcal{F}}: [0,1)^{\mathrm{op}} \to \operatorname{Shv}(C)$ (resp. $(0,1]^{\mathrm{op}} \to \operatorname{Shv}(C)$) valued in the
subcategory of monomorphisms extends uniquely to a left-continuous (resp. right-continuous) parabolic
sheaf, provided there is a commutative diagram
\[
\begin{tikzcd}[column sep=large]
    {\displaystyle\varprojlim_{0 \leq t < 1} \pa{\mathcal{F}}[t]}
        \arrow[r, "\sim"] \arrow[d, hook]
        & {\pa{\mathcal{F}}[0] \otimes \mathcal{O}_{C}(-D)} \arrow[d] \\
    {\pa{\mathcal{F}}[0]} \arrow[r, "\sim"]
        & {\pa{\mathcal{F}}[0] \otimes \mathcal{O}_{C}}
\end{tikzcd}
\]
in which the right vertical arrow is the canonical map (resp. the same diagram with $\pa{\mathcal{F}}[1]$ in place of
$\varprojlim_{0 \leq t < 1} \pa{\mathcal{F}}[t]$ and $\varinjlim_{0 < t \leq 1} \pa{\mathcal{F}}[t]$ in place of $\pa{\mathcal{F}}[0]$).
\end{remark}
\begin{definition}\label{def:dual-parabolic-sheaf}
    Let $\pa{\mathcal{F}}$ be a parabolic sheaf on $(C,D)$. The dual parabolic sheaf $\pa{\mathcal{F}}^{\vee}$ 
    is defined by 
    \[\pa{\mathcal{F}}^{\vee}[\alpha] := \Big(\varinjlim_{\beta > -\alpha} \pa{\mathcal{F}}[\beta]\Big)^{\vee}(-D)\]
    for $\alpha \in [0,1)$, with transition maps induced by the dual of the transition 
    maps of $\pa{\mathcal{F}}$.
\end{definition}

For completeness, we also define parabolic sheaves in the relative setting.

\begin{definition}\label{def:relative-parabolic-sheaf}
    Let $S$ be a $\mathbb{C}$-scheme of finite type, $c:\mathcal{C} \to S$ an $S$-curve 
    and $\mathcal{D}$ a relative Cartier divisor on $\mathcal{C}/S$. A parabolic sheaf 
    $\pa{\mathcal{F}}$ on $(\mathcal{C},\mathcal{D})/S$ is a functor 
    $\pa{\mathcal{F}}:\mathbb{R}^{\mathrm{op}} \to \operatorname{Shv}(\mathcal{C})$ such that
    $\pa{\mathcal{F}}[(-)+1] \simeq \pa{\mathcal{F}}[(-)] \otimes \mathcal{O}_{\mathcal{C}}(-\mathcal{D})$.
\end{definition}

\begin{remark}
    While all the results of this section are stated for parabolic sheaves
    on curves over $\mathbb{C}$, everything localizes as expected  \cite[Definition 1.0]{Yokogawa1993} \cite[Theorem 6.1]{BorneVistoli2012}. Later on, we will
    apply the results of this section to the relative case without making note of it.
\end{remark}

While our definitions of parabolic bundles and parabolic sheaves are standard, our definition of coparabolic bundles is not. We rectify this issue by proving that
 $\coparvect$ admits a fully faithful functor to $\pa{\operatorname{Shv}}(C,D)$ sending a coparabolic bundle $\widehat{\pa{V}}$ to the coparabolic bundle associated to $\pa{V}$ in the sense of \cite[p.~140]{Yokogawa1995}.

\begin{proposition}
    The mappings
    \begin{align*}
    \iota : \operatorname{Ob}(\parvect)\times [0,1) \to \operatorname{Vect}(C) &:
    (\pa{V}, \alpha) \mapsto \iota(\pa{V})[\alpha] := \ker \Big(V \to \bigoplus_{x_{j} \in D} V_{x_{j}}/V^{\min\{i \mid \alpha^{i}_{j} \geq \alpha\}}_{j} \Big) \\
    \widehat{\iota} : \operatorname{Ob}(\coparvect)\times (0,1] \to \operatorname{Vect}(C) &:
    (\widehat{\pa{W}}, \alpha) \mapsto \widehat{\iota}(\widehat{\pa{W}})[\alpha] := \ker \Big(W \to \bigoplus_{x_{j} \in D} W_{x_{j}}/W^{\max\{k \mid \beta^{k}_{j} \leq \alpha\}}_{j} \Big)
    \end{align*}
    respectively extend to fully faithful functors $\iota : \parvect \to \pa{\operatorname{Shv}}(C,D)$ and
    $\widehat{\iota} : \coparvect \to \pa{\operatorname{Shv}}(C,D)$, which act on morphisms by restriction,
    $\iota(f)[\alpha] = f|_{\iota(\pa{V})[\alpha]}$.
\end{proposition}

\begin{proof}
    If we let $\pa{Q_{V}}[\alpha] := \bigoplus_{x_{j} \in D} V_{x_{j}}/V^{\min\{i \mid \alpha^{i}_{j} \geq \alpha\}}_{j}$,
    we get a parabolic sheaf by declaring $\pa{Q_{V}}[\alpha\to \beta]$ to be the natural projection.
    It is well defined as the parabolic filtration is decreasing as parabolic weights increase, and then a simple diagram chase shows that $\iota(\pa{V})[\alpha \to \beta]$ exists and is uniquely determined by
    the values $\iota(\pa{V})[\alpha]$. 
    
    Likewise given parabolic bundles $\pa{V},\pa{V}'$ and a morphism $f: V \to V'$ 
    on the underlying vector bundles, by definition it is a parabolic morphism if and only if 
    it descends to $\pa{Q_{V}} \to \pa{Q_{V'}}$. Chasing the following commutative diagram 
    with exact rows shows that $f$ descends as above if and only if it restricts to a morphism 
    $\iota(\pa{V}) \to \iota(\pa{V}')$:
\[\begin{tikzcd}
	0 & {\iota(\pa{V})} & V & {\pa{Q_{V}}} & 0 \\
	0 & {\iota(\pa{V}')} & V' & {\pa{Q_{V'}}} & 0
	\arrow[from=1-1, to=1-2]
	\arrow[from=1-2, to=1-3]
	\arrow[dashed, from=1-2, to=2-2]
	\arrow[from=1-3, to=1-4]
	\arrow["f", from=1-3, to=2-3]
	\arrow[from=1-4, to=1-5]
	\arrow[dashed, from=1-4, to=2-4]
	\arrow[from=2-1, to=2-2]
	\arrow[from=2-2, to=2-3]
	\arrow[from=2-3, to=2-4]
	\arrow[from=2-4, to=2-5]
\end{tikzcd}\]
where the dashed arrows exist if and only if $f$ is a parabolic morphism.
Hence a fortiori both $\iota(-)$ and $Q_{(-)}$ are functorial, and $\iota(-)$ is fully faithful.
The same argument applies to the case of coparabolic bundles by replacing $\pa{Q_{V}}$ by $\widehat{\pa{Q_{V}}}:= \bigoplus_{x_{j} \in D} W_{x_{j}}/W^{\max\{k \mid \beta^{k}_{j} \leq \alpha\}}_{j}$ and putting hats on the parabolic objects throughout.
\end{proof}
\begin{proposition}\label{prop:coparabolic-ff}
    For a parabolic sheaf $\pa{\mathcal{F}}$, write $\widehat{\pa{\mathcal{F}}}[\alpha] := \varinjlim_{\beta > \alpha} \pa{\mathcal{F}}[\beta]$.
    The functors $\iota$ and $\widehat{\iota}$ are compatible with the equivalence of
    Proposition~\ref{prop:parabolic-coparabolic-equivalence} and with parabolic dualization, in the sense that there are
    natural isomorphisms
    \[
        \widehat{\iota(-)} \simeq \widehat{\iota}(\widehat{(-)})
        \qquad\text{and}\qquad
        \iota(-)^{\vee} \simeq \iota\big((-)^{\vee}\big)
    \]
    of functors $\parvect \to \pa{\operatorname{Shv}}(C,D)$.
\end{proposition}

\begin{proof}
    Both sides are literally the same functor, as
    Definitions~\ref{def:corresponding-coparabolic} and~\ref{def:parabolic-dual} are read
    off from $\iota(\pa{V})$. Passing to $\widehat{\iota(\pa{V})}$ replaces $\geq$ by $>$
    in the definition of $\iota$, so $\widehat{\iota(\pa{V})}[0] = V'$ and
    $\widehat{\iota(\pa{V})}[\alpha] \subseteq V'$ is cut out by the subspaces
    $\rho_{j}^{-1}(V^{\min\{i \mid \alpha^{i}_{j} > \alpha\}}_{j})$, which are the steps
    $V'^{\max\{k \mid \beta^{k}_{j} \leq \alpha\}}_{j}$ of $\widehat{\pa{V}}$ by
    construction. Likewise, dualizing turns kernels into kernels of annihilators, so the
    parabolic dual is the pointwise dual
    $\widehat{\iota}(\widehat{\pa{V}^{\vee}})[t] = \iota(\pa{V})[-t]^{\vee}(-D)$, and it
    turns the colimit of Definition~\ref{def:dual-parabolic-sheaf} into a limit. Hence
    \[
        \iota(\pa{V})^{\vee}[\alpha]
        = \varprojlim_{t < \alpha} \widehat{\iota}(\widehat{\pa{V}^{\vee}})[t]
        = \iota(\pa{V}^{\vee})[\alpha]
    \]
    by the first isomorphism.
\end{proof}

\begin{remark}
    \spaceskip=\fontdimen2\font plus \fontdimen3\font minus 1.4\fontdimen4\font
    In particular ${V^{\vee}[\alpha] \simeq (\widehat{V}[-\alpha])^{\vee}(-D)}$. This agrees with the conventions of \cite[Lemma 2.6.2]{landesman2023geometric} and \cite[Proposition 3.7]{Yokogawa1995}.
\end{remark}
Unless this causes confusion, we will suppress $\iota$ (resp. $\widehat{\iota}$) from our notation and use $\pa{V}$ (resp. $\widehat{\pa{V}}$)
to denote both a parabolic bundle (resp. a coparabolic bundle) and its associated parabolic sheaf $\iota(\pa{V})$. 

The category of parabolic sheaves has much better categorical properties. Grothendieck's axioms AB$5$ and AB$3^{*}$ are conditions ensuring good homological behavior. Recall that AB$5$ requires that the category has arbitrary coproducts with exact filtered colimits, while AB$3^{*}$ requires arbitrary products \cite[\S1.5]{Grothendieck1957}.
\begin{proposition}
    $\pa{\operatorname{Shv}}(C,D)$ is an abelian category which satisfies Grothendieck's axioms AB$5$ and AB$3^{*}$.
\end{proposition}

\begin{proof}
    The functor category $\operatorname{Shv}(C)^{\mathbb{R}^{\mathrm{op}}}$ is abelian and satisfies AB$5$ and AB$3^{*}$
    since $\operatorname{Shv}(C)$ does (see \cite[Proposition 3.1.1]{Grothendieck1957} or
    \cite[Tag 077P]{stacks-project}). Since $\mathcal{O}_{C}(-D)$ is flat and finitely presented, 
    it preserves all pointwise limits, pointwise colimits, and exactness of short exact sequences. 
    Therefore the inclusion $\pa{\operatorname{Shv}}(C,D) \hookrightarrow \operatorname{Shv}(C)^{\mathbb{R}^{\mathrm{op}}}$ 
    is exact and creates all limits and colimits in $\pa{\operatorname{Shv}}(C,D)$, hence the latter inherits AB$5$ and AB$3^{*}$.
\end{proof}

\begin{proposition}\label{parabolic-tensor-inner-hom}
    The category $\pa{\operatorname{Shv}}(C,D)$ of parabolic sheaves has a closed symmetric 
    monoidal structure. That is, for any parabolic sheaves $\pa{\mathcal{F}}, \pa{\mathcal{F}}' \in \operatorname{Ob}(\pa{\operatorname{Shv}}(C,D))$, 
    there exist a tensor product $\pa{\mathcal{F}} \otimes \pa{\mathcal{F}}'$ and an internal hom 
    $\pa{\ihom}(\mathcal{F}, \mathcal{F}')$ in $\pa{\operatorname{Shv}}(C,D)$, satisfying the usual 
    adjunction ${\Hom(\pa{\mathcal{G}} \otimes \pa{\mathcal{F}}, \pa{\mathcal{F}}')} \simeq {\Hom(\pa{\mathcal{G}}, \pa{\ihom}(\mathcal{F}, \mathcal{F}'))}$ 
    for all $\pa{\mathcal{G}} \in \operatorname{Ob}(\pa{\operatorname{Shv}}(C,D))$.
\end{proposition}
    
\begin{proof}
    We start with the tensor product.
    Given $\pa{\mathcal{F}},\pa{\mathcal{F}}' \in \operatorname{Ob}(\pa{\operatorname{Shv}}(C,D))$, the tensor product $\pa{\mathcal{F}}\otimes \pa{\mathcal{F}}'$ 
    is the parabolic sheaf defined by
    \[
    (\pa{\mathcal{F}}\otimes \pa{\mathcal{F}}')[\alpha] := \operatorname{coeq}\Big(
        \varinjlim_{\alpha \leq \beta + \beta' + 1} \pa{\mathcal{F}}[\beta] \otimes_{\mathcal{O}_{C}} \mathcal{O}_{C}(-D) \otimes_{\mathcal{O}_{C}} \pa{\mathcal{F}}'[\beta']
        \rightrightarrows
        \varinjlim_{\alpha \leq \beta + \beta'} \pa{\mathcal{F}}[\beta] \otimes_{\mathcal{O}_{C}} \pa{\mathcal{F}}'[\beta']
    \Big)
    \]
    
    In the category of presheaves, colimits are computed pointwise, and each 
    $\pa{\mathcal{F}}[\beta] \otimes_{\mathcal{O}_{C}} \pa{\mathcal{F}}'[\beta']$ is a sheaf. Since sheafification
    is a left adjoint, it preserves colimits, so this construction yields a sheaf for each $\alpha$.
    
    It remains to verify that $(\pa{\mathcal{F}}\otimes \pa{\mathcal{F}}')[\alpha+1] \simeq (\pa{\mathcal{F}} \otimes \pa{\mathcal{F}}')[\alpha]\otimes \mathcal{O}_{C}(-D)$. 
    But this follows from the fact that $\pa{\mathcal{F}}[\beta+1] \simeq \pa{\mathcal{F}}[\beta]\otimes \mathcal{O}_{C}(-D)$
    and $\pa{\mathcal{F}}'[\beta'+1] \simeq \pa{\mathcal{F}}'[\beta']\otimes \mathcal{O}_{C}(-D)$, together with the 
    compatibility of tensor products with colimits. The unit is the trivially parabolic $(\mathcal{O}_{C})_{\mathrm{triv}}$, 
    and the associativity and symmetry constraints are induced by those of $\otimes_{\mathcal{O}_{C}}$, since colimits commute with colimits.

    As to the matter of inner homs for parabolic bundles, if we define $\pa{\ihom}(\mathcal{F},\mathcal{F}') := \pa{\mathcal{F}}^{\vee} \otimes \pa{\mathcal{F}}'$, then
    \begin{align*}
       \pa{\ihom}(\mathcal{F},\mathcal{F}')[\alpha] &= \pa{\mathcal{F}}[-\alpha]^{\vee}\otimes \mathcal{O}_{C}(-D) \otimes \pa{\mathcal{F}}'[\alpha]\\
       \pa{\ihom}(-, \pa{\ihom}(\mathcal{F},\mathcal{F}'))[\alpha] &= (-)^{\vee}[\alpha] \otimes \pa{\mathcal{F}}[1-\alpha]\otimes \mathcal{O}_{C}(-D) \otimes \pa{\mathcal{F}}'[\alpha] = \pa{\ihom}(- \otimes \mathcal{F}, \mathcal{F}')[\alpha]
    \end{align*}
    as desired. For the general case, see \cite[\S3]{Yokogawa1995}.
\end{proof}

\begin{remark}
    If $\pa{V}$ and $\pa{V}'$ are both parabolic vector bundles, then we can explicitly 
    write their tensor product as
    \[
    (\pa{V}\otimes\pa{V}')[\alpha] =
    \sum_{\alpha \leq \beta + \beta' < \alpha + 1}\pa{V}[\beta] \otimes \pa{V}'[\beta']
    \subset (\pa{V} \otimes \pa{V}')[\alpha](D)
    \]
    and the sum has finitely many distinct summands. From this, we infer that 
    $\pa{V}\otimes\pa{V}'$ is itself a parabolic vector bundle.
\end{remark}

\begin{proposition}\label{prop:parabolic-serre-duality}
    Let $\pa{V}$ be a parabolic vector bundle on $(C,D)$. 
    Then for all $i \geq 0$, we have a canonical isomorphism
    \[H^{i}(C, \pa{V}) \simeq H^{1-i}(C, \widehat{\pa{V}^{\vee}} \otimes \omega_{C}(D))^{\vee}\]
\end{proposition}
\begin{proof}
    See \cite[Proposition 3.7]{Yokogawa1995}. For $\pa{V}$ with the trivial parabolic structure, $\widehat{\pa{V}^{\vee}}[0] = V^{\vee}(-D)$ and the statement is classical Serre duality.
\end{proof}
\begin{definition}
    Let $\pa{V}$ be a parabolic vector bundle. Its parabolic degree is 
    \[\pardeg \pa{V}:= \deg V + \sum^{n}_{j=1} \sum^{\len(\pa{V},j)}_{i=1}
    \alpha^{i}_{j} \dim(V^{i}_{j}/V^{i+1}_{j})\]
    and its parabolic slope is
    \[ \pmu(\pa{V}) := \frac{\pardeg \pa{V}}{\rk \pa{V}}\]
\end{definition}
\begin{definition}
    Given a parabolic bundle $\pa{V}$, $\deg V$ (resp. $\mu(V)$) denotes the degree (resp. the slope) of the underlying vector bundle.
\end{definition}
In the proposition below, we collect some basic facts about parabolic slopes that will be useful in later sections and which we use freely.
\begin{proposition}\label{basic-facts}
    Let $\pa{V},\pa{V}'$ be parabolic bundles of positive rank. Then 
    \begin{itemize}
        \item $\pmu(\pa{V}\otimes\pa{V}') = \pmu(\pa{V}) + \pmu(\pa{V}')$
        \item $\pmu(\pa{V}^{\vee}) = -\pmu(\pa{V})$
        \item $\pmu(\pa{\ihom}(V,V')) = \pmu(\pa{V}') - \pmu(\pa{V})$
        \item $\pardeg \pa{V} = \pardeg \pa{V}' + \pardeg \pa{V}''$ for every parabolic subbundle $\pa{V}' \subset \pa{V}$ with quotient $\pa{V}'' := \pa{V}/\pa{V'}$ carrying the induced parabolic structures.
    \end{itemize}
\end{proposition}
\begin{proof}
    These are straightforward computations from the definitions. See \cite[Lemma 2.6.5]{landesman2023geometric} 
    for the first two formulas. The third follows from the first two and the fact that 
    $\pa{\ihom}(V,V') \simeq V^{\vee} \otimes V'$ as parabolic bundles. For the last, see \cite[Part 3, \S I.A, Remark 3]{Seshadri1982}.
\end{proof}

\begin{lemma}\label{lemma:degree-set}
    Let $r \geq 0$ be an integer and, for each $j$, let $A_{j} \subset [0,1)$ be a finite set. The parabolic degrees of the parabolic quotient bundles of parabolic bundles of rank $r$ on $(C,D)$ whose parabolic weights at $x_{j}$ lie in $A_{j}$ form a closed discrete subset of $\mathbb{R}$.
\end{lemma}
\begin{proof}
    Let $\pa{Q}$ be such a quotient. Its induced weights $\beta^{i}_{j}$ at $x_{j}$ are among the weights of the bundle it is a quotient of, hence lie in $A_{j}$, and their multiplicities $\dim(Q^{i}_{j}/Q^{i+1}_{j})$ add up to $\rk Q \leq r$. Thus $\pardeg \pa{Q}$ lies in $\mathbb{Z} + \sum_{j=1}^{n} \sum_{\alpha \in A_{j}} c_{j,\alpha}\, \alpha$ for some integers $0 \leq c_{j,\alpha} \leq r$, and there are finitely many such translates of $\mathbb{Z}$.
\end{proof}

\begin{definition}
    Let $\pa{V}$ be a parabolic vector bundle. We say it is semistable (resp. stable) if 
    for every nonzero proper parabolic subbundle $\pa{V}' \subset \pa{V}$ the inequality 
    $\pmu(\pa{V}) \geq \pmu(\pa{V}')$ (resp. $\pmu(\pa{V}) > \pmu(\pa{V}')$) holds.
    We say that a coparabolic bundle is semistable (resp. stable) if its corresponding parabolic bundle is semistable (resp. stable).
\end{definition}

\begin{remark}
    Parabolic vector bundles with fixed rank, degree, and weight type form an algebraic stack. 
    Moreover, there exists a quasi-projective coarse moduli space $M^{\mathrm{par}}_{ss}(r,d,\alpha)$ 
    parametrizing S-equivalence classes of semistable parabolic bundles of rank $r$, degree $d$, 
    and weight type $\alpha = (\alpha^{i}_{j})$. When the weights are rational, this moduli space 
    can be constructed via geometric invariant theory. See \cite{MehtaSeshadri1980} for the 
    original construction in the case of parabolic bundles, and
    \cite[\S2, \S4]{Yokogawa1993} and
    \cite[\S5]{Yokogawa1995} or \cite{BodenYokogawa1996} for the general case. 
    The stack of all parabolic bundles (not necessarily semistable) is an Artin stack, 
    often denoted $\operatorname{ParBun}_{r,d}(C,D,\alpha)$.
\end{remark}

\section{The Atiyah bundle of a vector bundle}\label{sec:atiyah-bundle}

The Atiyah bundle parametrizes first-order differential operators on a vector bundle \cite[Theorem 1]{Atiyah1957}. Following \cite[\S3]{landesman2023geometric}, we recall the construction of Atiyah bundles for filtered parabolic bundles.

\begin{notation}
Throughout this section, let $C$ denote a smooth projective curve over $\mathbb{C}$ of genus $g \geq 0$, and let $D = x_{1} + \cdots + x_{n}$ be a reduced effective divisor on $C$. All sheaves are coherent $\mathcal{O}_{C}$-modules unless otherwise specified.

\begin{itemize}
\item $T_{C}$ denotes the tangent sheaf of $C$, the locally free sheaf of rank 1 dual to the sheaf of Kähler differentials $\Omega^{1}_{C}$.

\item $T_{C}(-D) := T_{C} \otimes \mathcal{O}_{C}(-D)$ denotes the subsheaf of $T_{C}$ consisting of vector fields vanishing along $D$.

\item $\Omega^{1}_{C}(\log D)$ denotes the sheaf of logarithmic 1-forms, characterized as the $\mathcal{O}_{C}$-module generated locally by $\Omega^{1}_{C}$ and $\{dz/z\}$ where $z$ is a local equation for $D$.

\item A filtration $F_{\bullet}: 0 = F_{0} \subset F_{1} \subset \cdots \subset F_{m} = V$ by subbundles means each $F_{i}$ is a locally free subsheaf with constant rank, and the inclusions are locally split. When the filtration is on a parabolic vector bundle, we endow each step of the filtration with the induced parabolic structure of Definition~\ref{def:induced-parabolic}.

\item $\operatorname{Art}_{/\mathbb{C}}$ denotes the category of local Artin $\mathbb{C}$-algebras with residue field $\mathbb{C}$.

\end{itemize}
\end{notation}

\subsection{Classical Atiyah bundles}
\begin{definition}
Let $V$ be a vector bundle over $C$.
Given an open set $U \subset C$, a section $f \in \mathcal{O}_{C}(U)$ and a $\mathbb{C}$-linear map $\tau \in \ihom_{\mathbb{C}}(V,V)(U)$, we define
\[
[f,\tau] : V(U) \to V(U) : s \mapsto \tau(f\cdot s) - f \cdot \tau(s).
\]
The sheaf of first-order differential operators is the sheafification of the presheaf
\[U \mapsto  \{ \tau \in \ihom_{\mathbb{C}}(V,V)(U) \mid \forall f \in \mathcal{O}_{C}(U),\; [f,\tau] \in \ihom_{\mathcal{O}_{C}}(V,V)(U) \}  \]
and is denoted $\operatorname{Diff}^{1}(V,V)$ \cite[\S16.8.1]{Gro65}.
\end{definition}
\begin{definition}
Let $V$ be a vector bundle over $C$. The Atiyah bundle $\aAt{V}$ is the subsheaf of $\operatorname{Diff}^{1}(V,V)$ whose sections over an open $U \subseteq C$ consist of those $\tau \in \operatorname{Diff}^{1}(V,V)(U)$ satisfying the following condition: for each $f \in \mathcal{O}_{C}(U)$, there exists a unique element $\delta_{\tau}(f) \in \mathcal{O}_{C}(U)$ such that
\[\tau(fs) - f\tau(s) = \delta_{\tau}(f) \cdot s\]
for all $s \in V(U)$, where the right-hand side denotes scalar multiplication of the section $s$ by $\delta_{\tau}(f)$.
\end{definition}

\begin{remark}
For $\tau \in \aAt{V}(U)$, the assignment $f \mapsto \delta_{\tau}(f)$ defines a $\mathbb{C}$-linear derivation $\delta_{\tau} : \mathcal{O}_{C}(U) \to \mathcal{O}_{C}(U)$. Indeed, for $f, g \in \mathcal{O}_{C}(U)$ and $s \in V(U)$:
\begin{align*}
\delta_{\tau}(fg) \cdot s &= \tau(fgs) - fg\tau(s)\\
&= \tau(fgs) - f\tau(gs) + f\tau(gs) - fg\tau(s)\\
&= (f\delta_{\tau}(g) + g\delta_{\tau}(f)) \cdot s
\end{align*}
This yields a morphism of coherent sheaves $\delta : \aAt{V} \to T_{C}$ whose kernel consists precisely of the $\mathcal{O}_{C}$-linear operators, hence we obtain the Atiyah exact sequence:
\begin{equation}\label{eq:atiyah-exact-sequence}
0 \to \End(V) \to \aAt{V} \xrightarrow{\delta} T_{C} \to 0
\end{equation}
\end{remark}

\begin{definition}
	Let $(V,F_{\bullet})$ be a filtered vector bundle on $C$. We respectively call the coherent subsheaves of filtration-preserving endomorphisms
	\[ \End(V,F_{\bullet}) := U \to \{ f \in \End(V)(U) | f(F_i) \subset F_i \text{ for all } 0 \leq i \leq m \}\]
	\[ \aAt{V,F_{\bullet}} := U \to \{ f \in \aAt{V}(U) | f(F_i) \subset F_i \text{ for all } 0 \leq i \leq m \}\]
	the endomorphism bundle and the Atiyah bundle of the filtered bundle $(V,F_{\bullet})$.
\end{definition}

\begin{remark}
It is immediate from the definitions that the subsheaf of $\mathcal{O}_{C}$-linear operators in $\aAt{V,F_{\bullet}}$ is just $\End(V,F_{\bullet})$, hence by restricting $\delta$ we obtain the filtered version of the Atiyah exact sequence:
\begin{equation}\label{eq:atiyah-filtered}
0 \to \End(V,F_{\bullet}) \to \aAt{V,F_{\bullet}} \xrightarrow{\delta} T_{C} \to 0
\end{equation}
\end{remark}

\begin{definition}
Let $(V,F_{\bullet})$ be a filtered vector bundle on $C$. The relative Atiyah bundle with respect to the divisor $D$ is the preimage
\[\At{V,F} := \delta^{-1}(T_{C}(-D)) \subset \aAt{V,F_{\bullet}}\]
where $T_{C}(-D) \subset T_{C}$ is the subsheaf of vector fields vanishing along $D$.
\end{definition}

\begin{remark}
By construction, $\At{V,F}$ fits in the short exact sequence
\begin{equation}\label{eq:relative-atiyah}
0 \to \End(V,F_{\bullet}) \to \At{V,F} \xrightarrow{\delta} T_{C}(-D) \to 0
\end{equation}
When $F_{\bullet}$ is the trivial filtration $(0 \subset V)$, we write $\At{V}$ for simplicity.
\end{remark}

\subsection{The parabolic Atiyah bundle}

To define the parabolic Atiyah bundle, recall that the sheaf $\pa{\End}(V) \subset \End(V)$ of Definition~\ref{def:parabolic-internal-hom} consists of endomorphisms preserving the parabolic flags at each $x_{j} \in D$ (see also \cite[Definition 1.5]{MehtaSeshadri1980} and \cite[Part 3, \S I.A, D\'efinition 3]{Seshadri1982}).

\begin{definition}\label{def:residue-map}
Let $\pa{V}$ be a parabolic vector bundle over $(C,D)$.

Since $\pa{\End}(V) \subset \End(V)$ is a subsheaf of locally free sheaves (Proposition~\ref{prop:parabolic-internal-hom-bundle}), the quotient sheaf $\End(V)/\pa{\End}(V)$ is coherent and supported on $D$. For each $d \in D$, the stalk $(\End(V)/\pa{\End}(V))_{d}$ is a finite-dimensional vector space. For each $d \in D$, define the residue map
\[\overline{r}_{d} : \At{V,F}|_{d} \to (\End(V)/\pa{\End}(V))_{d}\] 
as the composition 
\[\At{V,F}|_{d} \xrightarrow{r_{d}} \End(V,F_{\bullet})|_{d} \to \End(V)|_{d} \to (\End(V)/\pa{\End}(V))_{d},\]
where the last two maps are the natural inclusion and quotient \cite[\S2.2]{BiswasHeuHurtubise2020}, and $r_{d}$ is the residue homomorphism of \cite[\S3.2.3]{landesman2023geometric}. To define it, consider the commutative diagram with exact rows
\[\begin{tikzcd}[column sep=small]
0 \arrow[r] & \End(V,F_{\bullet})|_{d} \arrow[r] \arrow[d, equal] & \At{V,F}|_{d} \arrow[r, "\delta"] \arrow[d] & T_{C}(-D)|_{d} \arrow[r] \arrow[d, "0"] & 0 \\
0 \arrow[r] & \End(V,F_{\bullet})|_{d} \arrow[r] & \aAt{V,F_{\bullet}}|_{d} \arrow[r, "\delta"] & T_{C}|_{d} \arrow[r] & 0
\end{tikzcd}\]
obtained by restricting to $d$ the inclusion of \eqref{eq:relative-atiyah} into \eqref{eq:atiyah-filtered}. The right vertical map is the evaluation of the inclusion $T_{C}(-D) \hookrightarrow T_{C}$ at $d$, which is zero, so the middle one factors through $\End(V,F_{\bullet})|_{d}$ and $r_{d}$ is this factorization. Note that it restricts to the identity on $\End(V,F_{\bullet})|_{d}$.
\end{definition}

\begin{definition}\label{def:parabolic-atiyah}
The parabolic Atiyah bundle of $(\pa{V},F_{\bullet})$ is the kernel
\[\pAt{V,F} := \ker\left(\At{V,F} \xrightarrow{\oplus_{d \in D} \overline{r}_{d}} \bigoplus_{d \in D} (\End(V)/\pa{\End}(V))_{d}\right).\]
\end{definition}

The filtered parabolic Atiyah bundle fits in the short exact sequence
\begin{equation}\label{eq:atiyah-parabolic}
0 \to \pa{\End}(V,F_{\bullet}) \to \pAt{V,F} \xrightarrow{\delta} T_{C}(-D) \to 0
\end{equation}

In the case where the filtration $F_{\bullet}$ is trivial we suppress it from our notation, in which case the parabolic Atiyah exact sequence reads as
\begin{equation}\label{eq:parabolic-endomorphism}
0 \to \pa{\End}(V) \to \pAt{V} \xrightarrow{\delta} T_{C}(-D) \to 0
\end{equation}

\subsubsection{Deformation theory}

\begin{definition}
Let $(S,s_{0})$ be a pointed $\mathbb{C}$-scheme. A flat deformation of $(C,D)$ with base $(S,s_{0})$ is a cartesian square
\[\begin{tikzcd}[row sep=large, column sep=large]
(C,D) \arrow[r, hook] \arrow[d] \arrow[dr, phantom, "\lrcorner", very near start] & (\mathcal{C},\mathcal{D}) \arrow[d, "c"] \\
\{s_{0}\} \arrow[r, hook] & S
\end{tikzcd}\]
in which $c: \mathcal{C} \to S$ is flat of relative dimension $1$ and $\mathcal{D} \subset \mathcal{C}$ is a relative Cartier divisor flat over $S$.
A flat deformation of $(C,D,V)$ (resp. $(C,D,\pa{V})$) consists of a flat deformation of $(C,D)$ together with a vector bundle $\mathcal{V}$ on $\mathcal{C}$ (resp. parabolic bundle $\pa{\mathcal{V}}$ on $(\mathcal{C},\mathcal{D})$ relative to $S$, in the sense of Definition~\ref{def:relative-parabolic-sheaf}) whose restriction to the central fiber is isomorphic to $V$ (resp. $\pa{V}$).
\end{definition}

\begin{definition}\label{def:deformation-filtered}
A filtered parabolic bundle $(\pa{V},F_{\bullet})$ on $(C,D)$ is a parabolic bundle $\pa{V}$ together with a filtration $F_{\bullet}$ of $V$ by parabolic subbundles, each carrying the induced parabolic structure.
A flat deformation of $(C,D,V,F_{\bullet})$ (resp. $(C,D,\pa{V},F_{\bullet})$) consists of a flat deformation $\mathcal{V}$ of $(C,D,V)$ (resp. $\pa{\mathcal{V}}$ of $(C,D,\pa{V})$) together with a filtration $\mathcal{F}_{\bullet}$ of $\mathcal{V}$ by subbundles restricting to $F_{\bullet}$ on the central fiber (resp. and such that each $\mathcal{F}_{k}$ endowed with the parabolic structure induced from $\pa{\mathcal{V}}$ is a parabolic bundle on $(\mathcal{C},\mathcal{D})$ relative to $S$) \cite[Definition 3.5.4]{landesman2023geometric}.
\end{definition}

\begin{definition}
The deformation functor \cite[\S2.4]{Sernesi2006}
\[\operatorname{Def}_{(C,D,\pa{V},F_{\bullet})} : \operatorname{Art}_{/\mathbb{C}} \to \operatorname{Set}\] 
sends a local Artin $\mathbb{C}$-algebra $A$ to the set of isomorphism classes of flat deformations of $(C,D,\pa{V},F_{\bullet})$ over $\operatorname{Spec}(A)$. We define $\operatorname{Def}_{(C,D)}$ and $\operatorname{Def}_{(C,D,V,F_{\bullet})}$ in the same way, and omit $F_{\bullet}$ from the notation when it is trivial.
\end{definition}

\begin{lemma}\label{lem:locally-trivial}
Let $A = \mathbb{C}[\varepsilon]/\varepsilon^{2}$, write $U[\varepsilon] := U \times \operatorname{Spec} A$ for an open $U \subseteq C$, and let $\mathcal{U} = \{U_{i}\}_{i \in I}$ be an affine open cover of $C$ trivializing $V$, with $U_{ij} := U_{i} \cap U_{j}$ and transition functions $g_{ij}$ and $m_{ij}$ for $V$ and $\mathcal{O}_{C}$. Every first-order deformation of $(C,D,\pa{V},F_{\bullet})$ is isomorphic over each $U_{i}[\varepsilon]$ to $(U_{i}, D \cap U_{i}, \pa{V}|_{U_{i}}, F_{\bullet}|_{U_{i}}) \times \operatorname{Spec} A$ by an isomorphism restricting to the identity on the central fiber. In particular, it is determined up to isomorphism by transition functions
\[G_{ij} = g_{ij} + \varepsilon\tau_{ij} : U_{ij} \to \GL_{r}(A), \qquad f_{ij} = m_{ij} + \varepsilon\rho_{ij} : U_{ij} \to A^{*}\]
such that $G_{ij}$ preserves both $F_{\bullet} \otimes_{\mathbb{C}} A$ and the flags $V^{\bullet}_{d} \otimes_{\mathbb{C}} A$.
\end{lemma}

\begin{proof}
Denote the kernel of the homomorphism $\mathrm{GL}(V)(A) \to \mathrm{GL}(\mathbb{C})$ by $K(V)$. In terms of local trivializations over the affine opens $U_{i}$, a first-order deformation of $(C,D,\pa{V},F_{\bullet})$ is given by a tuple
\[(\{f_{ij}\}, \{G_{ij}\}, \{\mathrm{id}+\varepsilon\sigma_{i}\},\{\mathrm{id}+\varepsilon\eta_{i}\}) \in
\check{C}^{1}(\mathcal{U}, T_C(-D) \oplus \mathrm{GL}(V)(A)) \oplus \check{C}^0(\mathcal{U},K(V)^{\oplus 2}) 
\]
where each cocycle respectively controls the deformation of the pointed curve, the vector bundle, the filtration, and the parabolic structure. We write $G_{ij} =: g_{ij}+\varepsilon \tau_{ij}$ for convenience. 

We first straighten the deformation of the filtration, which is $(\mathrm{id} + \varepsilon\eta_{i})(F_{\bullet} \otimes_{\mathbb{C}} A)$ for some $\eta_{i} \in \End(V)(U_{i})$. The maps $\mathrm{id}+\varepsilon \eta_i : V \oplus \varepsilon V \to V \oplus \varepsilon V$ are invertible and so induce an isomorphism of deformations that we now make explicit: $\{G_{ij}\}$ transforms by conjugation and expanding we find that $\tau_{ij}$ gets replaced by $\tau_{ij} + g_{ij}\eta_{i} - \eta_{j}g_{ij}$, while the $0$-cocycles $\mathrm{id}+\varepsilon\eta_i$ transform by precomposition into $(\mathrm{id}+\varepsilon \eta_i)(\mathrm{id}-\varepsilon \eta_i) = \mathrm{id}$. 

By Definition~\ref{def:residue-map} the residue at $d$ changes by $(g_{ij}\eta_{i} - \eta_{j}g_{ij})|_{d}$ and the flag at $d$ becomes $(\mathrm{id} - \varepsilon\eta_{i})(V^{\bullet}_{d} \otimes_{\mathbb{C}} A)$. Both are harmless as soon as each $\eta_{i}|_{d}$ stabilizes $V^{\bullet}_{d}$, and since $\End(V,F_{\bullet})(U_{i})$ surjects onto $\bigoplus_{d \in D \cap U_{i}} \End(V,F_{\bullet})|_{d}$, whether we can choose such $\eta_{i}$ depends only on whether the class of $\eta_{i}$ vanishes in
\[\End(V)|_{d}/(\End(V,F_{\bullet}) + \pa{\End}(V))|_{d}\]
which we recognize to be the cokernel of the residue map $\overline{r}_{d}$.

The vanishing of this class is exactly the condition of Definition~\ref{def:deformation-filtered}. To see this, fix $d$, restrict everything to $d$, and write $\eta$ for $\eta_{i}|_{d}$. The $\ell$-th step of the flag $(\mathrm{id} - \varepsilon\eta)(V^{\bullet}_{d} \otimes_{\mathbb{C}} A)$ consists of the elements $u + \varepsilon(v - \eta u)$ with $u, v \in V^{\ell}_{d}$, and such an element lies in $F_{k} \otimes_{\mathbb{C}} A$ if and only if $u \in F_{k}$ and $v - \eta u \in F_{k}$. Their intersection therefore reduces onto $\{u \in F_{k} \cap V^{\ell}_{d} : \eta u \in F_{k} + V^{\ell}_{d}\}$ with kernel $\varepsilon(F_{k} \cap V^{\ell}_{d})$, so it is a direct summand lifting $F_{k} \cap V^{\ell}_{d}$ if and only if $\eta(F_{k} \cap V^{\ell}_{d}) \subset F_{k} + V^{\ell}_{d}$. A splitting of $V|_{d}$ adapted to both $F_{\bullet}$ and $V^{\bullet}_{d}$ shows that this inclusion holds for all $k$ and $\ell$ precisely when $\eta \in (\End(V,F_{\bullet}) + \pa{\End}(V))|_{d}$.
\end{proof}

\begin{proposition}\label{prop:deformation-atiyah}
Let $(V,F_{\bullet})$ (resp. $(\pa{V},F_{\bullet})$) be a filtered vector bundle (resp. filtered parabolic bundle, Definition~\ref{def:deformation-filtered}) on $C$ (resp. $(C,D)$). Then first-order deformations of $(C,D,\pa{V},F_{\bullet})$ (resp. $(C,D,V,F_{\bullet})$) are classified by $H^{1}(C,\pAt{V,F})$ (resp. $H^{1}(C,\aAt{V,F_{\bullet}})$).
\end{proposition}

\begin{proof}
For the reader's convenience, we include a proof of this well-known result. The correspondence between first-order deformations and Atiyah cohomology is classical; see \cite[\S3.3]{Sernesi2006} for the simple case of a (non-parabolic) line bundle. The parabolic generalization can be found in \cite[\S3.1]{BiswasHeuHurtubise2020}. To the author's knowledge no written-down explicit proof exists for the case of filtered parabolic bundles.

We only need to prove the parabolic case, as the non-parabolic case then follows by taking the trivial parabolic structure (Definition~\ref{def:trivial-parabolic-structure}). Fix $A = \mathbb{C}[\varepsilon]/\varepsilon^{2}$.

Let $\mathcal{U} = \{U_{i}\}_{i \in I}$ be an affine open cover of $C$ trivializing $V$, with $U_{ij} := U_{i} \cap U_{j}$ and $U_{ijk} := U_{i} \cap U_{j} \cap U_{k}$, and write $U[\varepsilon] := U \times \operatorname{Spec} A$ for an open $U \subseteq C$. Let $g_{ij}$ and $m_{ij}$ denote the transition functions for $V$ and $\mathcal{O}_{C}$ respectively. 
By Lemma~\ref{lem:locally-trivial}, every first-order deformation of $(C,D,\pa{V},F_{\bullet})$ is trivial over each $U_{i}[\varepsilon]$, compatibly with the filtration and the flags.
A first-order deformation of $(C,D,\pa{V},F_{\bullet})$ over $\operatorname{Spec} A$ consists of transition functions
\[G_{ij} = g_{ij} + \varepsilon\tau_{ij} : U_{ij} \to \GL_{r}(A), \qquad f_{ij} = m_{ij} + \varepsilon\rho_{ij} : U_{ij} \to A^{*}\]
lifting $g_{ij}$ and $m_{ij}$. 

The compatibility $G_{ij} \circ \iota_{i} = \iota_{j} \circ f^{*}_{ij}$, where $\iota_{i}: V|_{U_{i}} \to \mathcal{V}|_{U_{i}[\varepsilon]}$ denotes the trivialization, implies and is implied by
\[\tau_{ij}(fs) - m_{ij}(f)\tau_{ij}(s) = \rho_{ij}(f)g_{ij}(s) \qquad (f \in \mathcal{O}_{C}(U_{i}),\ s \in V(U_{i})),\]
which is equivalent to saying $\tau_{ij} \in \aAt{V}(U_{ij})$ with symbol $\delta(\tau_{ij}) = \rho_{ij}/m_{ij}$. Unwinding definitions shows that the gluing preserves the filtration and the flags if and only if
\[G_{ij}(F_{k} \otimes_{\mathbb{C}} A) = F_{k} \otimes_{\mathbb{C}} A, \qquad G_{ij}(\mathcal{V}^{k}_{d}) = \mathcal{V}^{k}_{d} \quad (d \in D \cap U_{ij}),\]
which is equivalent to $\tau_{ij} \in \pAt{V,F}(U_{ij})$.

The gluing carries the ideal $\mathcal{O}_{C}(-D) \otimes_{\mathbb{C}} A$ of $D[\varepsilon]$ into itself if and only if $\rho_{ij}(\mathcal{O}_{C}(-D)) \subset \mathcal{O}_{C}(-D)$, that is, if and only if $\delta(\tau_{ij}) \in T_{C}(-D)(U_{ij})$, which means $\tau_{ij} \in \At{V}(U_{ij})$.

Finally, as $\varepsilon^{2} = 0$, the cocycle relations $G_{ik} = G_{jk}G_{ij}$ and $f_{ik} = f_{ij}f_{jk}$ on $U_{ijk}$ amount to the 1-cocycle relation in $\check{C}^{1}(\mathcal{U}, \pAt{V,F})$,
\[\tau_{ik} = \tau_{jk}g_{ij} + g_{jk}\tau_{ij},\]
and applying $\delta$ gives the relation for $f_{ij}$. Read backwards, the same chain shows that any cocycle $\{\tau_{ij}\}$ defines a deformation through $G_{ij} := g_{ij} + \varepsilon\tau_{ij}$ and $f_{ij} := m_{ij}(1 + \varepsilon\delta(\tau_{ij}))$.

To establish injectivity, suppose two deformations yield cocycles $\{\tau_{ij}\}$ and $\{\tau'_{ij}\}$ differing by the coboundary of some $\{s_{i}\} \in \check{C}^{0}(\mathcal{U}, \pAt{V,F})$ so that $\tau'_{ij} - \tau_{ij} = s_{j}g_{ij} - g_{ij}s_{i}$. The isomorphism of deformations is given by $\Phi|_{U_{i}[\varepsilon]} = \mathrm{id} + \varepsilon s_{i}: \mathcal{V}|_{U_{i}[\varepsilon]} \to \mathcal{V}'|_{U_{i}[\varepsilon]}$, since on overlaps
\begin{align*}
(G'_{ij})^{-1} \circ \Phi|_{U_{j}} \circ G_{ij} &= (g_{ij} + \varepsilon\tau'_{ij})^{-1} \circ (\mathrm{id} + \varepsilon s_{j}) \circ (g_{ij} + \varepsilon\tau_{ij})\\
&= (\mathrm{id} - \varepsilon g^{-1}_{ij}\tau'_{ij})g^{-1}_{ij}(\mathrm{id} + \varepsilon s_{j})(g_{ij} + \varepsilon\tau_{ij})\\
&= \mathrm{id} + \varepsilon(- g^{-1}_{ij}\tau'_{ij} + g^{-1}_{ij}s_{j}g_{ij} + g^{-1}_{ij}\tau_{ij})\\
&= \mathrm{id} + \varepsilon s_{i} = \Phi|_{U_{i}}.
\end{align*}
Thus $\Phi$ glues to an isomorphism of deformations, which establishes the required bijection.
\end{proof}

\begin{corollary}\label{cor:filtration-obstruction}
Let $(\pa{V},F_{\bullet})$ be a filtered parabolic bundle over $(C,D)$, and $q \in H^{1}(C,\pAt{V})$ a first-order deformation of $(C,D,\pa{V})$. The obstruction to extending the filtration $F_{\bullet}$ along $q$ is the image $\pi_{F*}(q) \in H^{1}(C,\pa{\End}(V)/\pa{\End}(V,F_{\bullet}))$, and when it vanishes the extensions of $F_{\bullet}$ are parametrized by the quotient
\[H^{0}(C,\pa{\End}(V)/\pa{\End}(V,F_{\bullet}))/\im(H^{0}(C,\pa{\End}(V))).\]
\end{corollary}

\begin{proof}
From the exact sequence defining the parabolic Atiyah bundle with filtration, we obtain:
\begin{equation}\label{eq:filtration-exact}
0 \to \pAt{V,F} \xrightarrow{\iota} \pAt{V} \xrightarrow{\pi_{F}} \pa{\End}(V)/\pa{\End}(V,F_{\bullet}) \to 0
\end{equation}
The associated long exact sequence in cohomology gives:
\[\cdots \to H^{1}(C,\pAt{V,F}) \xrightarrow{\iota_{*}} H^{1}(C,\pAt{V}) \xrightarrow{\pi_{F*}} H^{1}(C,\pa{\End}(V)/\pa{\End}(V,F_{\bullet})) \to \cdots\]

By Proposition \ref{prop:deformation-atiyah}, the class $q \in H^{1}(C,\pAt{V})$ corresponds to a deformation of $(C,D,\pa{V})$. This deformation extends the filtration $F_{\bullet}$ if and only if $q$ lies in the image of $\iota_{*}$, which by exactness occurs precisely when $\pi_{F*}(q) = 0$. When the obstruction vanishes, the preimage $\iota^{-1}_{*}(q)$ is an affine space under the image of the connecting homomorphism%
\[\delta: H^{0}(C,\pa{\End}(V)/\pa{\End}(V,F_{\bullet})) \to H^{1}(C,\pAt{V,F}).\]

Since this map factors through the quotient by global endomorphisms preserving the filtration, we obtain the stated parametrization.
\end{proof}

\subsubsection{Connections and splittings}

\begin{definition}
Let $V$ be a vector bundle on $C$. A logarithmic connection is a $\mathbb{C}$-linear map $\nabla : V \to V \otimes \Omega^{1}_{C}(\log D)$ satisfying the Leibniz rule
\[\nabla(fs) = f\nabla(s) + s \otimes df\]
for local sections $f \in \mathcal{O}_{C}$ and $s \in V$.
\end{definition}

\begin{definition}
Let $V$ be a vector bundle on $C \setminus D$ and let $\nabla$ be a flat connection on $V$. We say the connection $\nabla$ has regular singularities if there exists an extension to a logarithmic connection $(\bar{V},\bar{\nabla})$ on $C$ with $\bar{\nabla}: \bar{V} \to \bar{V} \otimes \Omega^{1}_{C}(\log D)$ \cite[II, \S5.2]{Deligne1970}. We will refer to the pair $(\bar{V},\bar{\nabla})$ as a regular flat vector bundle.
\end{definition}

\begin{definition}
	Let $\pa{V}$ be a parabolic vector bundle over $(C,D)$. A regular flat connection on $\pa{V}$ is a section $q^{\nabla} : T_{C}(-D) \to \pAt{V}$ of the exact sequence \eqref{eq:parabolic-endomorphism}. We will refer to the pair $(\pa{V},\nabla)$ as a regular flat parabolic bundle.
\end{definition}
\begin{theorem}\label{thm:connection-parabolic-compatibility}
Let $(V,\nabla)$ be a logarithmic connection with regular singularities along $D$. For each $j \in \{1,\ldots,n\}$, let $\eta^{1}_{j}, \ldots, \eta^{s_{j}}_{j}$ denote the eigenvalues of the residue endomorphism $\operatorname{Res}_{x_{j}}(\nabla) \in \operatorname{End}(V|_{x_{j}})$, and set
\[
\alpha^{i}_{j} := \operatorname{Re}(\eta^{i}_{j}) - \lfloor\operatorname{Re}(\eta^{i}_{j})\rfloor \in [0,1).
\]
Reorder the $\alpha^{i}_{j}$ and remove repetitions so that $0 \leq \alpha^{1}_{j} < \alpha^{2}_{j} < \cdots < \alpha^{n_{j}}_{j} < 1$. Define
\[ V^{i}_{j} := \bigoplus_{\{\eta \,:\, \operatorname{Re}(\eta) - \lfloor\operatorname{Re}(\eta)\rfloor \,\geq\, \alpha^{i}_{j}\}} \ker\bigl((\operatorname{Res}_{x_{j}}(\nabla) - \eta)^{\rk V}\bigr). \]
Then the following hold:
\begin{enumerate}
\item $(V, \{V^{i}_{j}\}, \{\alpha^{i}_{j}\})$ defines a parabolic bundle $\pa{V}$.
\item The connection $\nabla$ induces a regular flat connection on $\pa{V}$.
\item $(V,\nabla)$ is the Deligne canonical extension of $(V,\nabla)|_{C \setminus D}$ if and only if $\operatorname{Re}(\eta) \in [0,1)$ for all eigenvalues $\eta$ of all $\operatorname{Res}_{x_{j}}(\nabla)$.
\end{enumerate}
\end{theorem}
\begin{proof}
See \cite[\S3.3]{landesman2023geometric} for details.
\end{proof}

\begin{proposition}\label{prop:connection-splitting}
Let $V$ be a vector bundle on $C$. There is a canonical bijection between
logarithmic flat connections $\nabla : V \to V \otimes \Omega^{1}_{C}(\log D)$ with regular singularities and splittings $q^{\nabla} : T_{C}(-D) \to \At{V}$ of the exact sequence \eqref{eq:relative-atiyah} \cite{BiswasHeuHurtubise2016}. The correspondence sends $\nabla$ to the map $q^{\nabla}$ defined by $q^{\nabla}(v)(s) = \nabla_{v}(s)$ for vector fields $v \in T_{C}(-D)$ and sections $s \in V$.
\end{proposition}

\begin{proof}
As $\Omega^{1}_{C}(\log D)$ is the dual of $T_{C}(-D)$, a $\mathbb{C}$-linear map $\nabla : V \to V \otimes \Omega^{1}_{C}(\log D)$ amounts to a $\mathbb{C}$-bilinear pairing $T_{C}(-D) \times V \to V$ which is $\mathcal{O}_{C}$-linear in the first variable. If we make $\ihom_{\mathbb{C}}(V,V)$ an $\mathcal{O}_{C}$-module through its target, adjunction turns such a pairing into an $\mathcal{O}_{C}$-linear map
\[q^{\nabla} : T_{C}(-D) \to \ihom_{\mathbb{C}}(V,V), \qquad q^{\nabla}(v)(s) = \nabla_{v}(s).\]
We claim $\nabla \mapsto q^{\nabla}$ is the bijection.

Adjunction being itself bijective, we only need to match the conditions on the two sides. Evaluated on a vector field $v$, the Leibniz rule reads
\[q^{\nabla}(v)(fs) - f\,q^{\nabla}(v)(s) = v(f)\,s,\]
which says precisely that $q^{\nabla}(v)$ lies in $\At{V}$ with symbol $v$. Thus $\nabla$ is a logarithmic connection if and only if $q^{\nabla}$ splits \eqref{eq:relative-atiyah}. Flatness and regularity come for free, since $\Omega^{2}_{C} = 0$ and a logarithmic connection is regular by definition.
\end{proof}

\begin{remark}
By the correspondence in Proposition \ref{prop:connection-splitting}, a regular flat connection on a parabolic bundle $\pa{V}$ is equivalent to a logarithmic connection $\nabla: V \to V \otimes \Omega^{1}_{C}(\log D)$ on the underlying bundle such that the induced residue endomorphisms at points of $D$ preserve the parabolic flags.
\end{remark}

\begin{remark}\label{rem:associated-local-system}
Suppose $(V,\nabla)$ is the Deligne canonical extension of $(V,\nabla)|_{C \setminus D}$, and give $V$ the parabolic structure of Theorem~\ref{thm:connection-parabolic-compatibility}. We call the resulting $(\pa{V},\nabla)$ the regular flat parabolic bundle associated to the local system $\ker(\nabla|_{C \setminus D})$.
\end{remark}

\begin{lemma}\label{lem:iso-section}
Let $(\pa{V},\nabla)$ be a regular flat parabolic bundle over $(C,D)$. The connection induces a section 
\[\operatorname{iso}_{\nabla}: \operatorname{Def}_{(C,D)}(\mathbb{C}[\varepsilon]/\varepsilon^{2}) \to \operatorname{Def}_{(C,D,\pa{V},\nabla)}(\mathbb{C}[\varepsilon]/\varepsilon^{2})\]
of the forgetful morphism.
\end{lemma}

\begin{proof}
By Proposition \ref{prop:deformation-atiyah}, deformations of $(C,D)$ and $(C,D,\pa{V})$ are classified by $H^{1}(C,T_{C}(-D))$ and $H^{1}(C,\At{\pa{V}})$ respectively. The forgetful map corresponds to the morphism 
\[\delta_{*}: H^{1}(C,\pAt{V}) \to H^{1}(C,T_{C}(-D))\] induced by the symbol map. The connection $\nabla$ gives a splitting 
$q^{\nabla}: T_{C}(-D) \to \At{\pa{V}} $
inducing $(q^{\nabla})_{*}: H^{1}(C,T_{C}(-D)) \to H^{1}(C,\pAt{V})$ with $\delta_{*} \circ (q^{\nabla})_{*} = \mathrm{id}$, which provides the desired section $\operatorname{iso}_{\nabla}$. For a more detailed exposition, see \cite[Proposition 3.5.7]{landesman2023geometric}.
\end{proof}

\section{Isomonodromic deformations}\label{sec:isomonodromy}
We review the analytic notion of isomonodromic deformation and introduce formal and algebraic analogues, proving their first-order equivalence (Theorem~\ref{serenity}).
 
\subsection{Analytic isomonodromic deformations}
\begin{definition}
    Let $\mathscr{C}$ be a complex manifold, $\Delta$ a contractible complex manifold, $\pi : \mathscr{C} \to \Delta$ a proper holomorphic 
    submersion of relative dimension $1$, and $d_{1},\ldots, d_{n} : \Delta \to \mathscr{C}$ disjoint sections. 
    Set $\mathscr{D} := \bigcup^{n}_{i=1} \im(d_{i})$. 
    Given a basepoint $0 \in \Delta$ and a regular flat vector bundle $(V,\nabla)$ on $(C,D) := (\pi^{-1}(0),\mathscr{D} \cap \pi^{-1}(0))$, an isomonodromic deformation is an extension to a flat bundle $(\mathscr{V},\widetilde{\nabla})$ on $(\mathscr{C},\mathscr{D})$ with logarithmic singularities along $\mathscr{D}$, restricting to a Deligne canonical extension on each fiber.
\end{definition}

\begin{lemma}
    With the notation above, for every regular flat vector bundle $(V,\nabla)$ there exists a unique isomonodromic deformation up to unique isomorphism.
\end{lemma}
\begin{proof}
By Ehresmann's theorem \cite{ehresmann1952connexions}, the morphism $\pi: \mathscr{C} \to \Delta$ is a locally trivial $C^{\infty}$ fibration. 
Since $\Delta$ is contractible, there exists a diffeomorphism $\Phi: (\mathscr{C},\mathscr{D}) \to (C \times \Delta, D \times \Delta)$ 
covering the identity on $\Delta$.

The pullback $\Phi^{*}((V,\nabla) \boxtimes \mathcal{O}_{\Delta})$ defines a flat bundle on $\mathscr{C}$ and a flat connection with 
logarithmic singularities along $\mathscr{D}$. By the Riemann--Hilbert correspondence, flat bundles
with regular singularities on $(\mathscr{C}, \mathscr{D})$ correspond to representations of 
$\pi_{1}(\mathscr{C} \setminus \mathscr{D})$. Since $\Delta$ is contractible, the inclusion 
$C \setminus D \hookrightarrow \mathscr{C} \setminus \mathscr{D}$ induces an isomorphism on fundamental groups, so any two isomonodromic deformations can only differ by the choice of extension of $\nabla$ to all of $\mathscr{C}$. But by \cite[II, \S5.4]{Deligne1970}, the choice of extension is uniquely determined by the requirement that $\widetilde{\nabla}$ be the Deligne canonical extension on each fiber.
\end{proof}

\begin{notation}\label{moduli-notation}
Let $X := \mathcal{M}_{g,n}$ denote the moduli stack of $n$-pointed smooth projective curves of genus $g$.
Let $Y := \mathcal{M}_{dR}(C,D,r)$ denote the de Rham moduli stack over $X$ of rank $r$ logarithmic connections
on curves of type $(C,D)$ with regular singularities, and $\varrho: Y \to X$ the forgetful morphism.
Denote the universal cover of $X$ by $\tilde{X}$. 
Let $\mathcal{M}_{B}(C,D,r)$ denote the Betti moduli stack over $X$ parametrizing rank $r$ representations
of $\pi_{1}(C \setminus D)$ up to conjugation, or equivalently rank $r$ local systems on $C \setminus D$. The Riemann--Hilbert correspondence induces
an analytic map $\mathcal{M}_{dR}(C,D,r)^{an} \to \mathcal{M}_{B}(C,D,r)^{an}$. Even though this map is not injective when $D$ is not empty, it admits local sections given by sending a representation of $\pi_{1}(C \setminus D)$ to the Deligne canonical extension of its corresponding connection on $C \setminus D$.
\end{notation}

\begin{definition}
   In the case where $\Delta := \tilde{X}$ is the universal cover of the moduli stack $\mathcal{M}^{an}_{g,n}$ and $\mathscr{C}$ is the universal curve over it, we call $(\mathscr{V},\widetilde{\nabla})$ 
   the universal analytic isomonodromic deformation, and refer to the image of the induced morphism $\mathscr{L}_{V}: \Delta \to Y/\Delta : s \mapsto [(\mathscr{V},\widetilde{\nabla})_{\mathscr{C}_{s}}]$ as the analytic isomonodromic leaf of $(V,\nabla)$.
\end{definition}

\begin{notation}
	Let $Z$ be a scheme or an analytic space. Given an admissible sheaf of ideals $\mathscr{I} \subset \mathcal{O}_{Z}$, 
	we denote the formal completion of $Z$ along $\mathscr{I}$ by $\operatorname{Spf}(Z,\mathscr{I})$. Given $z \in Z$, 
	we abbreviate $\operatorname{Spf}(Z,\mathfrak{m}_{z})$ by $\operatorname{Spf}(z)$.
\end{notation}

\begin{definition}\label{def:quot}
	Let $S$ be either a scheme over $Y$ or an analytic space over $Y^{an}$, write $(C_{S}, D_{S}, \pa{V_{S}})$ for the pulled-back family, and fix a polynomial $P$. The Quot functor on $S$-schemes, respectively analytic spaces over $S$, is
	\[ \Quot^{P}_{S}(T) := \left\{ \pa{V_{T}} \twoheadrightarrow \pa{Q} \;\middle|\; \begin{array}{l} \pa{Q} \text{ a quotient parabolic sheaf, flat over } T, \text{ with} \\ \text{Hilbert polynomial } P \text{ on the fibres with respect to } \omega_{C_{T}/T}(D_{T}) \end{array} \right\} \Big/ {\simeq}. \]
\end{definition}

A quotient of $\pa{V_{T}}$ carries the induced parabolic structure, so it is the same thing as a quotient of $V_{T}$, and the parabolic structure plays no role in the following theorem. When $S$ is a noetherian scheme the curve $C_{S} \to S$ is projective, so $\Quot^{P}_{S}$ is represented by a projective $S$-scheme by Grothendieck's theorem \cite[Theorem~5.14]{Nitsure2005}.

\begin{theorem}\label{thm:quot}
	For every analytic space $S$ over $Y^{an}$, the functor $\Quot^{P}_{S}$ is represented by an analytic space proper over $S$. When $S = U^{an}$ for a scheme $U$ of finite type over $Y$, this space is $(\operatorname{Quot}^{P}_{U})^{an}$.
\end{theorem}
\begin{proof}
	Since $C_{S} \to S$ is proper, by Pourcin's theorem \cite[Th\'eor\`eme~2]{Pourcin1969} we represent the functor of $T$-flat quotients of $V_{S}$ by an analytic space over $S$. We can find $\Quot^{P}_{S}$ as a connected component of said analytic space, because Hilbert polynomials are locally constant in proper flat families. Analytification of flat families of quotients of $V_{T}$ preserves the conditions on the algebraic and analytic Quot functors, thus we get a morphism $(\operatorname{Quot}^{P}_{U})^{an} \to \Quot^{P}_{U^{an}}$ of analytic spaces over $U^{an}$.

	Let $A$ be an Artinian local $\mathbb{C}$-algebra and $\operatorname{Spec} A \to U$ a point. The curve $C_{A}$ is a projective scheme over $\mathbb{C}$, so by GAGA \cite[\S12, Th\'eor\`emes~2 et~3]{Serre1956} analytification is an exact equivalence between coherent sheaves on $C_{A}$ and on $C_{A}^{an}$, compatible with restriction to the closed fibre. It therefore matches the $A$-flat quotients with Hilbert polynomial $P$ on the two sides. With $A = \mathbb{C}$ this says our morphism is bijective, and in general it says that it induces isomorphisms of completed local rings, since these pro-represent the Artinian points lying below a given $\mathbb{C}$-point. But a morphism $B \to B'$ of analytic local $\mathbb{C}$-algebras which becomes an isomorphism after completion is already one. It is injective because $B \to \widehat{B}$ is, and $\mathfrak{m}_{B}B' = \mathfrak{m}_{B'}$ because completion is faithfully flat, so $B'$ is finite over $B$ by the Weierstrass preparation theorem, and Nakayama's lemma applied to $B' = B + \mathfrak{m}_{B}B'$ gives $B' = B$.
\end{proof}

\subsection{Formal isomonodromic deformations}
\begin{definition}
	Let $C$ be a smooth projective curve over $\mathbb{C}$ and $D \subset C$ a reduced divisor. Further let $(A,\mathfrak{m})$ be a complete local ring with residue field $\mathbb{C}$. A formal $(A,\mathfrak{m})$ deformation of the pair $(C,D)$ is a  pullback square
	\[\begin{tikzcd}
		{(C,D)} & {(\widehat{\mathscr{C}},\widehat{\mathscr{D}})} \\
		{\operatorname{Spec} \mathbb{C}} & {\mathrm{Spf}(A,\mathfrak{m})}
		\arrow[from=1-1, to=1-2]
		\arrow[from=1-1, to=2-1]
		\arrow["\lrcorner"{anchor=center, pos=0.125}, draw=none, from=1-1, to=2-2]
		\arrow[from=1-2, to=2-2]
		\arrow[from=2-1, to=2-2]
	\end{tikzcd}\]
	where the bottom horizontal arrow is induced from $A \to A/\mathfrak{m}$, and $\widehat{\mathscr{C}}$ and $\widehat{\mathscr{D}}$ are flat over $\mathrm{Spf}(A,\mathfrak{m})$. We will often leave $(A,\mathfrak{m})$ implicit and refer to a pair $(\widehat{\mathscr{C}},\widehat{\mathscr{D}})$ that fits in a diagram like the one above as a formal deformation of $(C,D)$.
\end{definition}
\begin{definition}
	Let $(\pa{V},\nabla)$ be a regular flat parabolic bundle on $(C,D)$. A formal $(A,\mathfrak{m})$ isomonodromic deformation of $(\pa{V},\nabla)$ is a quadruple consisting of a formal $(A,\mathfrak{m})$ deformation $(\widehat{\mathscr{C}}, \widehat{\mathscr{D}})$ of $(C,D)$ together with a formal parabolic bundle $\widehat{\pa{\mathscr{V}}}$ and an absolutely flat connection $\widehat{\nabla} : \widehat{\pa{\mathscr{V}}} \to \widehat{\pa{\mathscr{V}}} \otimes \Omega^{1}_{\widehat{\mathscr{C}}/\mathbb{C}}(\log \widehat{\mathscr{D}})$ such that their restriction to the central fiber is $(\pa{V},\nabla)$.
\end{definition}

\begin{proposition}[cf.\ {\cite[Proposition~(8.9)]{Katz1970}}]\label{prop:horizontal-sections}
	Let $k$ be a field of characteristic zero, $R$ a $k$-algebra, $B := R[[t_{1},\dots,t_{m}]]$ and $I := (t_{1},\dots,t_{m})$. Let $M$ be an $I$-adically complete $B$-module with a flat connection $\nabla \colon M \to M \otimes_{B} \widehat{\Omega}^{1}_{B/R}$ relative to $R$. Then restriction to $t = 0$ induces a bijection $M^{\nabla = 0} \to M/IM$, and $M \simeq M^{\nabla=0} \widehat{\otimes}_{R} B.$
\end{proposition}
\begin{proof}
	 For convenience we abbreviate $\nabla_{\partial_{t_{j}}}$ by $\nabla_j$. We proceed by induction. Let $s^{(N)} \in M$ and define $E_{j}^{(N)} := \nabla_{j}s^{(N)}$ and assume $E_{j}^{(N)} \in
	I^{N}M$. Flatness of the connection gives us
	\[ t_{j} \nabla_{k} E_{j} = t_{j} \nabla_{j} E_{k} + I^{N+1}M \]
	hence after setting $s^{(N+1)} := s^{(N)}- \frac{1}{N+1} \sum_{j} t_{j} E_{j}$ simply computing shows that
	\[ \nabla_{k} s^{(N+1)} \equiv E_{k} - \tfrac{1}{N+1}\Big(E_{k} + \sum_{j} t_{j}\nabla_{j}E_{k}\Big) \equiv 
	E_{k} - \frac{1}{N+1}(E_{k}+NE_{k}) = 0 \pmod{I^{N+1}} \] and it is easy to see it is the only possible choice for $s^{(N+1)}$. Therefore there exists a unique extension of $s = s^{(0)} \in M/IM$ to a horizontal section of $M$, establishing the bijection we sought.
	
	As for the rest, $M^{\nabla=0} \twoheadrightarrow M/IM$ and $I^{N}/I^{N+1} = Sym^{N}(I/I^{2})$ as $R$-modules, thus by induction and the 5-lemma the middle vertical arrow of
	\[\begin{tikzcd}
		0 & {M^{\nabla=0}\otimes_{R} I^{N}/I^{N+1}} & {M^{\nabla=0}\otimes_{R}B/I^{N+1}} & {M^{\nabla=0} \otimes_{R} B/I^{N}} & 0 \\
		0 & {I^{N}M/I^{N+1}M} & {M/I^{N+1}} & {M/I^{N}} & 0
		\arrow[from=1-1, to=1-2]
		\arrow[from=1-2, to=1-3]
		\arrow[from=1-2, to=2-2]
		\arrow[from=1-3, to=1-4]
		\arrow[from=1-3, to=2-3]
		\arrow[from=1-4, to=1-5]
		\arrow[from=1-4, to=2-4]
		\arrow[from=2-1, to=2-2]
		\arrow[from=2-2, to=2-3]
		\arrow[from=2-3, to=2-4]
		\arrow[from=2-4, to=2-5]
	\end{tikzcd}\]
	is an isomorphism, and since $M$ is $I$-complete, passing to the limit establishes the surjectivity of the inclusion $M^{\nabla=0} \widehat{\otimes}_R B \to M$.
\end{proof}
\begin{proposition}\label{formal-isomonodromic-local}
	Let $k$ be a field of characteristic zero, $R$ a $k$-algebra, $B := R[[t_{1},\dots,t_{m}]]$ and $I := (t_{1},\dots,t_{m})$. Restriction to $t = 0$ is an equivalence from the category of $I$-adically complete $B$-modules with an absolute flat connection $\nabla : M \to M \otimes_{B} \widehat{\Omega}^{1}_{B/k}$ to the category of $R$-modules with a flat connection $\bar{\nabla} : N \to N \otimes_{R} \Omega^{1}_{R/k}$, with quasi-inverse $(N,\bar{\nabla}) \mapsto (N \widehat{\otimes}_{R} B,\, \bar{\nabla} \otimes 1 + 1 \otimes d)$. The same holds for log connections with poles along a locally principal ideal $J \subset R$ and along $JB$. In particular two absolute flat connections $(M,\nabla)$ and $(M',\nabla')$ extending the same $(M/I,\bar{\nabla})$ are isomorphic, by a unique isomorphism inducing the identity on $M/I$.
\end{proposition}
\begin{proof}
	As before by the 5-lemma and induction, $M = (M/I)\widehat{\otimes}_{R} B$. In the resp. case, if $J \subset R$ is a locally principal ideal and the connection on $M/I$ is valued in $\Omega^{1}_{R/k}(\operatorname{log}(J))$, then $JB \simeq J\widehat{\otimes}_{R} B$ by the same argument as for $M$, so without loss of generality we may assume that after subtracting by $\nabla_{R/k}\widehat{\otimes}_{R} \operatorname{id}_{B} $ the two connections have no poles.
	
	Either way, $\nabla-\nabla_{R/k}$ and $\nabla'-\nabla_{R/k}$ are relative connections $M \to M \otimes \Omega^{1}_{B/R}$, therefore we can apply Proposition \ref{prop:horizontal-sections} to get an isomorphism $\xi: M^{\nabla=0} \otimes \operatorname{id}_{B} \simeq M^{'\nabla'=0} \otimes \operatorname{id}_{B}$ compatible with the absolute connections. Likewise a horizontal map $f : M \to M'$ preserves the sections killed by the relative connections, so it is $\bar{f} \,\widehat{\otimes}\, \operatorname{id}_{B}$ for its restriction $\bar{f}$ to $t = 0$, and every map of this form is horizontal.
\end{proof}
\begin{lemma}\label{formal-isomonodromic}
	Let $(\pa{V},\nabla)$ be a regular flat parabolic bundle on $(C,D)$, and $(\widehat{\mathscr{C}},\widehat{\mathscr{D}})$ a formal $(A,\mathfrak{m})$ deformation of $(C,D)$ with $A = \mathbb{C}[[t_{1},\dots,t_{m}]]$. Then $(\pa{V},\nabla)$ has a unique formal $(A,\mathfrak{m})$ isomonodromic deformation over $(\widehat{\mathscr{C}},\widehat{\mathscr{D}})$ up to a unique isomorphism restricting to the identity on the central fiber.
\end{lemma}
\begin{proof}
	On $C \setminus D$ the flat connection $\nabla$ is the same thing as a stratification of $V$ \cite[Propositions 2.9 and 2.11, Theorem 2.15]{BerthelotOgus1978}, that is, an isomorphism
	\[ \xi : p_{2}^{*}(V,\nabla) \simeq p_{1}^{*}(V,\nabla) \]
	over the formal neighborhood $\widehat{\Delta_{C \setminus D}}$ of the diagonal, restricting to the identity on the diagonal and to $\nabla$ on its first-order neighborhood, and satisfying the cocycle condition $p_{13}^{*}\xi = p_{12}^{*}\xi \circ p_{23}^{*}\xi$ over the formal neighborhood of the diagonal in the triple product. It is universal: if two maps $f$ and $g$ from a formal scheme to $C \setminus D$ agree on a subscheme of definition, then $(f,g)$ factors through $\widehat{\Delta_{C \setminus D}}$, and $(f,g)^{*}\xi$ identifies $g^{*}(V,\nabla)$ with $f^{*}(V,\nabla)$. We claim that it pulls back to the descent datum of the isomonodromic deformation we seek.

	Let $\{U_{\alpha}\}$ be an affine open cover of $C$ such that each point of $D$ lies in exactly one $U_{\alpha}$, and write $\widehat{U}_{\alpha}$ for the open formal subscheme of $\widehat{\mathscr{C}}$ supported on $U_{\alpha}$. Since $U_{\alpha}$ is smooth and affine, the deformation $(\widehat{U}_{\alpha},\widehat{\mathscr{D}})$ of $(U_{\alpha}, D \cap U_{\alpha})$ is trivial (compare Lemma \ref{lem:locally-trivial}), that is, there is a retraction $\rho_{\alpha} : \widehat{U}_{\alpha} \to U_{\alpha}$ with $\rho_{\alpha}^{-1}(D) = \widehat{\mathscr{D}}$. For $\alpha \neq \beta$ the overlap $U_{\alpha\beta}$ lies in $C \setminus D$, and there $\rho_{\alpha}$ and $\rho_{\beta}$ need not agree, yet they do agree on the central fiber. We therefore have a commutative diagram
	\[\begin{tikzcd}[column sep=large]
		{\bigsqcup_{\alpha \neq \beta} \widehat{U}_{\alpha\beta}} & {\widehat{\Delta_{C \setminus D}}} \\
		{\bigsqcup_{\alpha} \widehat{U}_{\alpha}} & C \\
		{\widehat{\mathscr{C}}}
		\arrow["{(\rho_{\alpha},\rho_{\beta})}", from=1-1, to=1-2]
		\arrow[shift right, from=1-1, to=2-1]
		\arrow[shift left, from=1-1, to=2-1]
		\arrow["{p_{1}}"', shift right, from=1-2, to=2-2]
		\arrow["{p_{2}}", shift left, from=1-2, to=2-2]
		\arrow["{\rho_{\alpha}}"', from=2-1, to=2-2]
		\arrow[from=2-1, to=3-1]
	\end{tikzcd}\]
	whose left column is the cover of $\widehat{\mathscr{C}}$ by the $\widehat{U}_{\alpha}$ together with its overlaps.

	Pulling $\xi$ back along the top arrow gives isomorphisms
	\[ \varphi_{\alpha\beta} := (\rho_{\alpha},\rho_{\beta})^{*}\xi : \rho_{\beta}^{*}(V,\nabla) \simeq \rho_{\alpha}^{*}(V,\nabla) \]
	over the $\widehat{U}_{\alpha\beta}$, restricting to the identity on the central fiber. The pullback of a cocycle is a cocycle, so they form a descent datum. Moreover $\rho_{\alpha}^{*}\pa{V}$ is a parabolic bundle on $(\widehat{U}_{\alpha},\widehat{\mathscr{D}})$ because $\rho_{\alpha}$ pulls $D$ back to $\widehat{\mathscr{D}}$, and the $\varphi_{\alpha\beta}$ are isomorphisms of parabolic bundles since the overlaps do not meet $\widehat{\mathscr{D}}$. The $\rho_{\alpha}^{*}(\pa{V},\nabla)$ thus descend to an isomonodromic deformation of $(\pa{V},\nabla)$ over $(\widehat{\mathscr{C}},\widehat{\mathscr{D}})$.

	It is the only one. Indeed $\rho_{\alpha}$ identifies $\mathcal{O}(\widehat{U}_{\alpha})$ with $\mathcal{O}(U_{\alpha})[[t_{1},\dots,t_{m}]]$, so by Proposition \ref{formal-isomonodromic-local} the restriction of any isomonodromic deformation to $\widehat{U}_{\alpha}$ is isomorphic to $\rho_{\alpha}^{*}(\pa{V},\nabla)$, by a unique isomorphism restricting to the identity on $U_{\alpha}$. Its descent datum then consists of isomorphisms $\rho_{\beta}^{*}(V,\nabla) \simeq \rho_{\alpha}^{*}(V,\nabla)$ restricting to the identity on $U_{\alpha\beta}$, of which the same proposition allows only one, namely $\varphi_{\alpha\beta}$. The local isomorphisms therefore glue to the unique isomorphism we sought.
\end{proof}

\begin{lemma}\label{formal-existence}
	Let $(\pa{V}, \nabla)$ be a regular flat parabolic bundle over $(C, D)$, with
	corresponding $x \in X$ and $V \in Y$. There is a formal section
	\[
	L^{\wedge} : \operatorname{Spf}\widehat{\mathcal{O}}_{X,x}
	\to \operatorname{Spf}\widehat{\mathcal{O}}_{Y,V}
	\]
	of $\varrho^{\wedge}$ whose tangent map at $V$ extends $q^{\nabla}$.
\end{lemma}
\begin{proof}
	Apply Lemma \ref{formal-isomonodromic} to the universal formal deformation $(\widehat{\mathscr{C}},\widehat{\mathscr{D}})$ of $(C,D)$ over $\operatorname{Spf}\widehat{\mathcal{O}}_{X,x}$, where $\widehat{\mathcal{O}}_{X,x} = \mathbb{C}[[s_{1},\dots,s_{l}]]$. Composing $\widehat{\nabla}$ with $\Omega^{1}_{\widehat{\mathscr{C}}/\mathbb{C}}(\log\widehat{\mathscr{D}}) \to \Omega^{1}_{\widehat{\mathscr{C}}/\widehat{X}}(\log\widehat{\mathscr{D}})$ turns $\widehat{\pa{\mathscr{V}}}$ into a formal family of regular flat parabolic bundles on $(\widehat{\mathscr{C}},\widehat{\mathscr{D}})$ restricting to $(\pa{V},\nabla)$ on the central fiber. By the moduli property of $Y$ this family is classified by a morphism $L^{\wedge} \colon \operatorname{Spf}\widehat{\mathcal{O}}_{X,x} \to \operatorname{Spf}\widehat{\mathcal{O}}_{Y,V}$ over $X$.

	The tangent map can be read off from the gluing modulo $(s)^{2}$. In the notation used in the proof right above we can write
		\[ \rho_{\beta}^{\#} \equiv \rho_{\alpha}^{\#} + \sum_{j} s_{j}\, \rho_{\alpha}^{\#} \circ v^{\alpha\beta}_{j} \pmod{(s)^{2}} \qquad v^{\alpha\beta}_{j} \in T_C(-D)(U_{\alpha\beta}), \]
	and the cocycle $(v^{\alpha\beta}_{j})$ represents the Kodaira--Spencer class of $\partial_{s_{j}}$ in $H^{1}(C,T_{C}(-D))$. Modulo $(s)^{2}$ the map $(\rho_{\alpha},\rho_{\beta})$ factors through the first-order neighborhood of the diagonal, where $\xi$ is $\nabla$, hence
		\[ \varphi_{\alpha\beta}(\rho_{\beta}^{*}e) \equiv \rho_{\alpha}^{*}e + \sum_{j} s_{j}\, \rho_{\alpha}^{*}\big(\nabla_{v^{\alpha\beta}_{j}}e\big) \pmod{(s)^{2}}. \]
	The first-order deformation along $\partial_{s_{j}}$ is therefore glued along $\mathrm{id} + s_{j}\,q^{\nabla}(v^{\alpha\beta}_{j})$ in the notation of Proposition \ref{prop:connection-splitting}, the image under $q^{\nabla}$ of the Kodaira--Spencer cocycle. Under the classification of Proposition \ref{prop:deformation-atiyah} this says that $L^{\wedge}_{*}$ is the section $\operatorname{iso}_{\nabla}$ of Lemma \ref{lem:iso-section} on Zariski tangent spaces.
\end{proof}

Let $\pi: \mathscr{C} \to X$ denote the universal curve with sections $\sigma_{i}: X \to \mathscr{C}$ marking $\mathscr{D} = \sum^{n}_{i=1} \sigma_{i}(X)$. For any point $x \in X$ corresponding to $(C,D)$, denote the associated analytic universal isomonodromic deformation of $V$ by $\mathscr{V}$. Then by flatness 
\[\widehat{\mathscr{V}} := \varprojlim_{k} \left( \mathscr{V} \otimes \mathcal{O}_{X,x}/\mathfrak{m}^{k}_{x} \right) 
= \mathscr{V} \otimes \varprojlim_{k} \mathcal{O}_{X,x}/\mathfrak{m}^{k}_{x}\] 
is a vector bundle over the formal completion $\widehat{\mathscr{C}} := \operatorname{Spf}(\mathscr{C}, \pi^{*}\mathfrak{m}_{x})\to \operatorname{Spf}(x)$. 
The construction is independent of the choice of lift of $x$ as the resulting bundles are isomorphic   
via the monodromy action of $\pi_{1}(\mathcal{M}_{g,n})$.

\begin{definition}
    We will refer to the formal parabolic bundle $\widehat{\pa{\mathscr{V}}}$ over $(\widehat{\mathscr{C}},\widehat{\mathscr{D}})$ 
    as the formal universal isomonodromic deformation of $(\pa{V},\nabla)$, and to the associated morphism
    $\widehat{\mathscr{L}}_{V} : \operatorname{Spf}(x) \to \operatorname{Spf}(Y,\mathfrak{m}_{V})$ as the formal isomonodromic leaf of $(\pa{V},\nabla)$.
\end{definition}

\begin{remark}\label{formal-analytic-comparison}
	Since $C$ is projective, formal GAGA \cite[III, \S5]{EGA3} provides an equivalence between the categories of coherent sheaves on $\widehat{\mathscr{C}}$ for the Zariski and complex-analytic topologies and we use the two viewpoints interchangeably.
\end{remark}

\subsection{First-order isomonodromic deformations}
\begin{remark}\label{kodaira-spencer-map}
The connection $\nabla$ induces the Kodaira--Spencer map $q^{\nabla} : T_{C}(-D) \to \pAt{V}$ and given some filtration $F_{\bullet}$ the induced map on cohomology $H^{1}(q^{\nabla})$ factors through $H^{1}(C,\pAt{V,F})$ if and only if the filtration extends to first-order isomonodromic deformations by Corollary \ref{cor:filtration-obstruction}. By Deligne's extension theorem \cite[II, \S5.4]{Deligne1970} or Lemma \ref{lem:iso-section}, an infinitesimal deformation of $(C,D)$ lifts to one of $(V,\nabla)$, thus giving a splitting $q^{\nabla}$ of $\varrho_{*}$.
\end{remark}

\begin{definition}
	Let $A := \mathcal{O}_{X,x}/\mathfrak{m}^{2}_{x} = \mathbb{C}\oplus H^{1}(C,T_{C}(-D))^{\vee}$. Then 
	\[\Hom_{\mathbb{C}}(H^{1}(C,T_{C}(-D)),T_{V}Y) = 
	\Hom_{\text{sch}}(\operatorname{Spec}(A), \operatorname{Spec}(\mathcal{O}_{Y}/\mathfrak{m}^{2}_{V}))\]
	and we define the universal first-order isomonodromic deformation of the regular flat parabolic bundle 
	$(\pa{V},\nabla)$ to be the parabolic bundle $(\pa{\mathcal{V}},\nabla)$ over $\mathcal{C} := \mathscr{C} \times_{X} \operatorname{Spec}(A)$ determined by the 
	morphism $q^{\nabla} : \operatorname{Spec}(A) \to \operatorname{Spec}(\mathcal{O}_{Y}/\mathfrak{m}^{2}_{V})$.
\end{definition}
\begin{proposition}\label{lift-alg-formal}
    Let $(\pa{V}, \nabla)$ be a regular flat parabolic bundle over $(C, D)$. The formal
    isomonodromic leaf $\widehat{\mathscr{L}}_{V}$ is the moduli-classifying morphism
    of the formal isomonodromic deformation of $(\pa{V},\nabla)$, and its first-order truncation
    is the universal first-order isomonodromic deformation of $(\pa{V}, \nabla)$.
\end{proposition}
\begin{proof}
    By definition, $\widehat{\mathscr{L}}_{V}$ is the moduli-classifying morphism of the formal completion $\widehat{\mathscr{V}}^{an}$ of the analytic universal isomonodromic deformation $\mathscr{V}$. Isomonodromy gives $\widehat{\mathscr{V}}^{an}$ an absolutely flat connection restricting to $\nabla$ on the central fiber, hence it is canonically isomorphic to the formal isomonodromic deformation of Lemma~\ref{formal-isomonodromic}. Its tangent map at $V$ is $q^{\nabla}$ by Lemma~\ref{formal-existence}, so its restriction to a first-order neighborhood agrees with the definition we gave for the first-order isomonodromic deformation, as it should.
\end{proof}

\begin{theorem}\label{serenity}
    Let $(\pa{V},\nabla)$ be a regular flat parabolic bundle over $(C,D)$ and let $\pa{V'} \subset \pa{V}$ be a parabolic subbundle. As before we let $x \in X$ be the point corresponding to $(C,D)$, and let $\pa{\mathscr{V}}$, $\pa{\widehat{\mathscr{V}}}$, $\pa{\mathcal{V}}$ denote the analytic, formal, and first-order isomonodromic deformations of $(\pa{V},\nabla)$, with $\pa{\mathcal{V}}$ over $\operatorname{Spec}(A)$.
    Then the following are equivalent:
    \begin{itemize}
        \item[1.] There exists a neighborhood $U \subset \mathcal{M}^{an}_{g,n}$ of $(C,D)$
        such that $\pa{V'}$ extends to a parabolic subbundle $\pa{\mathscr{V}'} \hookrightarrow \pa{\mathscr{V}}|_{U}.$
        \item[2.] There exists a parabolic bundle $\pa{\widehat{\mathscr{V}'}}$ over $\widehat{\mathscr{C}}$ such that 
        $\pa{\widehat{\mathscr{V}'}} \hookrightarrow \pa{\widehat{\mathscr{V}}}$ and 
        $\pa{\widehat{\mathscr{V}'}}\otimes_{\mathcal{O}_{X}} k(x) \simeq \pa{V'}.$
    \end{itemize}
    and both imply
    \begin{itemize}
        \item[3.] There exists a parabolic bundle $\pa{\mathcal{V}'}$ over $\mathcal{C}$ such that 
        $\pa{\mathcal{V}'} \hookrightarrow \pa{\mathcal{V}}$ and 
        $\pa{\mathcal{V}'} \otimes_{A} k(x) \simeq \pa{V'}.$
    \end{itemize}

\end{theorem}
\begin{proof}
    $1 \implies 2$: by taking the formal completion along $C$, and $2 \implies 3$ by Proposition \ref{lift-alg-formal}. 
    
    $2 \implies 1$: the analyticity of $\mathscr{V}$ reduces the problem to finitely many unknowns, after which Artin's approximation theorem applies. Fix a relatively ample line bundle $\mathcal{O}(1)$ on $\mathscr{C}$, for instance $\omega_{\mathscr{C}/X}(\mathscr{D})$, and choose $m$ such that $V(m)$ and $V'(m)$ have vanishing $H^{1}$ and $V'(m)$ is globally generated. After shrinking $U$, Grauert's theorems \cite[\S7, S\"atze~3, 5, 6]{Grauert1960} make $E := \pi_{*}\mathscr{V}(m)$ a holomorphic vector bundle on $U$ whose formation commutes with base change \cite[Bemerkung~4.4.1]{KiehlVerdier1971}, and the vanishing of $H^{1}(C,V'(m))$ makes $\widehat{\pi}_{*}\widehat{\mathscr{V}'}(m)$ a subbundle of its completion. Thus $\widehat{\mathscr{V}'}$ defines a formal section $\widehat{s}$ of the Grassmann bundle $G \to U$ of $h^{0}(V'(m))$-planes in $E$.

    Let $S \subset E_{G}$ be the tautological subbundle on $G$, and let $Q$ be the cokernel of the evaluation map $S \otimes \mathcal{O} \to \mathscr{V}(m)$ on $\mathscr{C}_{G}$. Let $Z$ (resp. $\pa{Z}$) be the subfunctor whose points are the planes of $G$ generating subbundles of rank $r' := \rk V'$ (resp. which satisfy the conditions of Definition \ref{def:deformation-filtered}). Since formation of Fitting ideals commute with base change, the $T$-points of $Z$ are the locus where $\operatorname{Fitt}_{r-r'}(Q_{T})$ is the unit ideal and $\operatorname{Fitt}_{r-r'-1}(Q_{T})$ vanishes. The former is an open condition on $G$, as $\mathscr{C}_{G} \to G$ is proper. The latter says that the $(r'+1)$-st exterior power of the evaluation map vanishes on $\mathscr{C}_T$. This corresponds to a section $\phi$ of a vector bundle $M$ on $\mathscr{C}_{G}$. Twisting by $\mathcal{O}(k)$ with $k$ large enough to ensure that $\mathcal{O}(k)$ is generated by $\pi_{*}\mathcal{O}(k)$ and $M(k)$ has no $H^{1}$ on the fibers, this amounts to the vanishing of the induced map of vector bundles $\pi_{*}\mathcal{O}(k) \to \pi_{*}M(k)$ on $G$, which is a closed condition. In the parabolic case, Definition~\ref{def:deformation-filtered} requires the cokernels of $S \oplus \mathscr{V}^{i}_{d}(m) \to \mathscr{V}(m)|_{d}$ to be locally free of constant rank, which is again a locally closed condition on Fitting ideals.

    Over every infinitesimal neighborhood of $x$, the plane $\widehat{s}$ generates $\widehat{\mathscr{V}'}(m)$ by Nakayama's lemma, so $\widehat{s}$ factors through $\pa{Z}$. Near $\widehat{s}(x)$ the subspace $\pa{Z}$ is cut out of an open subset of $U \times \mathbb{C}^{N}$ by finitely many convergent power series, and Artin's theorem \cite[Theorem 1.2]{Artin1968} provides, for every $c \geq 1$, a holomorphic section $s_{c}$ of $\pa{Z}$ over a neighborhood of $x$ which agrees with $\widehat{s}$ modulo $\mathfrak{m}^{c}_{x}$. Since $V'(m)$ is globally generated, the plane $s_{c}(x) = H^{0}(C,V'(m))$ generates $V'(m)$, so twisting the subbundle generated by $s_{c}$ by $\mathcal{O}(-m)$, we obtain the desired extension $\pa{\mathscr{V}'} \hookrightarrow \pa{\mathscr{V}}$ of $\pa{V'}$. 
\end{proof}

\section{An infinitesimal criterion for extending filtrations}

Throughout this section, every regular flat parabolic bundle is the one associated to a local system on $C \setminus D$, in the sense of Remark~\ref{rem:associated-local-system}.

\begin{lemma}
    \label{lemma:ggg-criterion}
    Let $\pa{V}$ be a parabolic vector bundle on a smooth projective curve $(C, D)$ of genus $g$, with $\alpha := \max_{i, j} \alpha^{i}_{j}$. Let $V' \subset V$ be a subbundle, $c := \rk V - \rk V'$, $\delta := h^{0}(V) - h^{0}(V')$, and suppose $0 \leq \mu(\gr^{HN}_{i}(V')) \leq 2g - 2$ for every Harder--Narasimhan graded piece of $V'$. Then for any $\mu \in \mathbb{Q}$ with $\pmu(\pa{V}) \geq \mu + n\alpha$,
    \[
    \rk V \cdot (2g - 1 - \mu) \geq cg - \delta.
    \]
\end{lemma}

\begin{proof}
    By \cite[Lemma 6.2.1]{landesman2023geometric} applied to $V'$ and Riemann--Roch on $V$,
    \[
    \frac{\deg V'}{2} + \rk V' \geq h^{0}(V') = h^{0}(V) - \delta \geq \deg V - (g-1)\rk V - \delta.
    \]
    The hypothesis $\mu(\gr^{HN}_{i}(V')) \leq 2g-2$ gives $\mu(V') \leq 2g-2$, hence $\deg V' \leq (2g-2) \rk V'$ and $\deg V'/2 + \rk V' \leq g \rk V'$, whence
    \[
    g \rk V' \geq \deg V - (g-1)\rk V - \delta.
    \]
    By \cite[Lemma 6.3.3]{landesman2023geometric}, $\mu(V) \geq \pmu(\pa{V}) - n\alpha \geq \mu$, so $\deg V \geq \mu \rk V$. Substituting $\rk V' = \rk V - c$:
    \[
    g(\rk V - c) \geq \mu \rk V - (g-1) \rk V - \delta,
    \]
    and rearranging gives $(2g - 1 - \mu) \rk V \geq cg - \delta$.
\end{proof}

It can be more convenient to work with a version of Lemma~\ref{lemma:ggg-criterion} with no hypothesis on the Harder--Narasimhan filtration of the subbundle. The following lemma is \cite[Lemma 4.4]{landesman2024big} with the slope equality relaxed to an inequality.

\begin{lemma}
    \label{lemma:non-ggg-parabolic}
    Let $\widehat{\pa{V}}$ be a semistable coparabolic vector bundle on $(C, D)$ and $\mu \in \mathbb{Q}$ with $\pmu(\widehat{\pa{V}}) \geq \mu + n$ and $\mu < 2g - 1$. Let $U \subset \widehat{\pa{V}}[0]$ be a proper subbundle, $c := \rk V - \rk U$, and $\delta := h^{0}(C, \widehat{\pa{V}}[0]) - h^{0}(C, U)$. Then
    \[
    \rk V \geq \frac{cg - \delta}{2g - 1 - \mu}.
    \]
\end{lemma}
\begin{proof}
    By \cite[Lemma 6.3.4]{landesman2023geometric} applied to $\widehat{\pa{V}}$, every quotient vector bundle of $\widehat{\pa{V}}[0]$ has slope at least $\pmu(\widehat{\pa{V}}) - n \geq \mu$. We conclude by applying \cite[Proposition 4.3]{landesman2024big} to the proper subbundle $U \subset \widehat{\pa{V}}[0]$.
\end{proof}

\begin{corollary}
    \label{cor:ggg-criterion-coparabolic}
    Let $\widehat{\pa{V}}$ be a coparabolic vector bundle on $(C, D)$, and let $V'$, $c$, $\delta$ satisfy the hypotheses of Lemma~\ref{lemma:ggg-criterion}. Then for any $\mu \in \mathbb{Q}$ with $\pmu(\widehat{\pa{V}}) \geq \mu + n$,
    \[
    \rk V \cdot (2g - 1 - \mu) \geq cg - \delta.
    \]
\end{corollary}
\begin{proof}
    The underlying parabolic bundle $\pa{V}$ has the same parabolic slope. Since $\alpha < 1$, $\pmu(\widehat{\pa{V}}) \geq \mu + n \geq \mu + n\alpha$, and Lemma~\ref{lemma:ggg-criterion} applies.
\end{proof}

\begin{definition}
Let $V$ be a vector bundle on a smooth curve $C$ and $\mathcal{F} \subset V$ a coherent subsheaf. The saturation $\widetilde{\mathcal{F}}$ of $\mathcal{F}$ in $V$ is the kernel of the composition
\[V \twoheadrightarrow V/\mathcal{F} \twoheadrightarrow (V/\mathcal{F})/\mathrm{tors}\]
where $(V/\mathcal{F})/\mathrm{tors}$ denotes the quotient by the torsion subsheaf.
\end{definition}

\begin{lemma}\label{lemma:filtration-connection}
Let $(\pa{V}, \nabla)$ be a regular flat parabolic bundle on $(C, D)$,
and let $0 = F_{0}\pa{V} \subset F_{1}\pa{V} \subset F_{2}\pa{V} = \pa{V}$ be a 2-step filtration by parabolic subbundles such that $\nabla(F_{1}\pa{V}) \not\subset F_{1}\pa{V} \otimes \omega_{C}(D)$, and write $\pa{G} := \pa{\ihom}(\gr^{F}_{1}(V),\gr^{F}_{2}(V))$. 
Let $q^{\nabla} : T_{C}(-D) \to \pa{\End}(V)/\pa{\End}(V,F_{\bullet})$ be the map induced by the connection, and assume $H^{1}(q^{\nabla}) = 0$. Then the
composite
\[T_{C}(-D) \xrightarrow{q^{\nabla}} \pa{\End}(V)/\pa{\End}(V,F_{\bullet}) \xrightarrow{\pi} \pa{G}\]
is nonzero and vanishes on $H^{1}$.
\end{lemma}
\begin{proof}
	The surjection $\pi: \pa{\End}(V)/\pa{\End}(V,F_{\bullet}) \twoheadrightarrow \pa{G}$ is the composition of restriction along $F_{1}\pa{V} \hookrightarrow \pa{V}$ in the source and projection
	along $\pa{V} \twoheadrightarrow \gr^{F}_{2}(\pa{V})$ in the target. Thus $H^{1}(\pi \circ q^{\nabla}) = H^{1}(\pi) \circ H^{1}(q^{\nabla}) = 0$, and that the composite $\pi \circ q^{\nabla}$ is nonzero is a restatement of the hypothesis that $\nabla$ does not preserve $F_{1}\pa{V}$.
\end{proof}
\begin{lemma}\label{lemma:quantitative-two-step}
	Let $(\pa{V}, \nabla)$ be a regular flat parabolic bundle on $(C, D)$, and let $0 = F_{0}\pa{V} \subset F_{1}\pa{V} \subset F_{2}\pa{V} = \pa{V}$ be a $2$-step filtration by parabolic subbundles such that $\gr^{F}_{i}(\pa{V})$ are semistable of positive rank and $\nabla(F_{1}\pa{V}) \not\subset F_{1}\pa{V} \otimes \omega_{C}(D)$. Further write $\pa{G} := \pa{\ihom}(\gr^{F}_{1}(V), \gr^{F}_{2}(V))$, $r_{i} := \rk(\gr^{F}_{i}(\pa{V}))$, and assume $d := \pardeg(\gr^{F}_{2}(\pa{V})) \geq 0$. If the composite
	\[
	q^{\nabla}: T_{C}(-D) \xrightarrow{q^{\nabla}} \pa{\End}(V)/\pa{\End}(V,F_{\bullet}) \xrightarrow{\pi} \pa{G}
	\]
	vanishes on $H^{1}$, then $r_{1}r_{2} + d(r_{1}+r_{2}) \geq g$.
	
\end{lemma}
\begin{proof}
	Since $\pardeg \pa{V} = 0$ \cite[B.3]{EsnaultViehweg1986}, additivity of parabolic degrees implies $\gr^{F}_{i}(\pa{V})$ have slopes $\frac{(-1)^{i} d}{r_{i}}$. Thus $\pa{G}$ is semistable \cite[Theorem 2 and Main Theorem]{Simpson1990} of parabolic slope
	\[
	s := \frac{d}{r_{1}} + \frac{d}{r_{2}} = \frac{d(r_{1}+r_{2})}{r_{1}r_{2}} \geq 0,
	\]
	and $\pa{G}^{\vee} \otimes \omega_{C}(D)$ is semistable of parabolic slope $2g-2+n-s$.
	
	By Lemma \ref{lemma:filtration-connection}, the composite $q^{\nabla} := \pi \circ q^{\nabla} \colon T_{C}(-D) \to \pa{G}$ is nonzero and vanishes on $H^{1}$. Its Serre dual, by Proposition \ref{prop:parabolic-serre-duality}, is a nonzero coparabolic morphism $\widehat{\pa{G}^{\vee}} \otimes \omega_{C}(D) \to \omega^{\otimes 2}_{C}(D)$ vanishing on global sections. The global sections of $E := \big(\widehat{\pa{G}^{\vee}} \otimes \omega_{C}(D)\big)[0]$ therefore generate a subsheaf whose saturation $U$ is a proper subbundle with $h^{0}(C, U) = h^{0}(C, E)$. Applying Lemma \ref{lemma:non-ggg-parabolic} with $\mu = 2g-2-s$, we obtain with $c \geq 1$ and $\delta = 0$ that
	\[
	r_{1}r_{2} = \rk \pa{G} \geq \frac{g}{2g-1-\mu} = \frac{g}{1+s}.
	\]
	Multiplying through by $1+s$ yields our desired inequality.
\end{proof}

\section{Main results}

Throughout this section, every regular flat parabolic bundle is the one associated to a local system on $C \setminus D$, in the sense of Remark~\ref{rem:associated-local-system}.

\begin{theorem}\label{thm:stability}
Let $(\pa{V},\nabla)$ be an irreducible regular flat parabolic bundle of rank $r$ on $(C,D)$, and use Lemma \ref{lemma:degree-set} to enumerate the degrees of nonnegative quotients of $\pa{V}$: $0 \leq d_0 < d_1 \cdots$ . If $\lfloor r^2/4 \rfloor + d_k r < g$ for some $k \in \mathbb{N}$, then an analytically general isomonodromic deformation of $\pa{V}$ has no nonzero proper parabolic subbundle of parabolic degree at least $-d_k$. 
Moreover, every proper parabolic subbundle $\pa{V'}$ of parabolic degree $-d_{k+1}$ of an analytically general isomonodromic deformation is stable, and satisfies $r_1 r_2 + d_{k+1} (r_1+r_2) \geq g$, where $r_1 := \rk \pa{V'}$ and $r_2 := r - r_1$.
\end{theorem}
\begin{proof}
By \cite[Corollary 6.1.2]{landesman2023geometric}, an analytically general isomonodromic deformation of $\pa{V}$ is semistable, so we may assume without loss of generality that $\pa{V}$ is parabolically semistable and $\{-d_k\}_{k=0}^{\infty}$ enumerates the parabolic degrees of its parabolic subbundles. 

Once we fix the rank, underlying degree, and flag intersection dimensions at $D$, the parabolic subbundles of the fibers of $\pa{\mathscr{V}}$ with these invariants form a locally closed subspace of $\Quot^{P}_{\Delta}$, where $P$ is the Hilbert polynomial of the quotient and $\Delta$ is viewed as an analytic space over $Y^{an}$ through the leaf $\mathscr{L}_{V}$, and its closure is proper over $\Delta$ by Theorem~\ref{thm:quot}. Since locally on $\Delta$ there are only finitely many choices of these discrete invariants of parabolic degree at least $-d_k$, if an analytically general isomonodromic deformation of $\pa{V}$ has such a subbundle, then one of these locally closed subspaces must dominate an open neighborhood $U \subset \Delta$ of $(C,D)$. After perhaps shrinking $U$, taking a local section yields a family
of parabolic subbundles $\pa{\mathscr{V}}' \subset \pa{\mathscr{V}}|_{U}$ of parabolic degree at least $-d_k$, which by Theorem~\ref{serenity} extends to first order in every isomonodromic direction.

So suppose after passing to another analytically general curve, every nonzero proper parabolic subbundle of parabolic degree at least $-d_k$ is semistable. Let $\pa{V}'$ be such a subbundle of rank $r_1$, write $r_2 := r - r_1$, and suppose it persists along the isomonodromic deformation near some $t\in U$. The filtration $0 \subset \pa{V}' \subset \pa{\mathscr{V}_t}$ extends to first order in every direction, so the composite $\pi \circ q^{\nabla}$ vanishes on $H^1$ by Corollary \ref{cor:filtration-obstruction}, and it is nonzero by the irreducibility of $\nabla$. Lemma \ref{lemma:quantitative-two-step} applies and gives 
	\[r_1 r_2 + d_k(r_1 + r_2) \geq g > \lfloor (r_1+r_2)^2/4 \rfloor + d_k(r_1+r_2) \geq r_1 r_2 + d_k(r_1 + r_2), \]
which is absurd. From this we deduce that an analytically general isomonodromic deformation has no nonzero proper parabolic subbundle of parabolic degree at least $-d_k$.

Without loss of generality we may pass to yet another analytically general curve $\pa{\mathscr{V}_t}$ and freely assume it has no nonzero proper parabolic subbundle of parabolic degree at least $-d_k$. We claim any subbundle $\pa{V}' \subset \pa{\mathscr{V}_t}$ of parabolic degree $-d_{k+1}$ is stable. A nonzero proper parabolic subbundle $\pa{V}^{''} \subset \pa{V}'$ is a parabolic subbundle of $\pa{\mathscr{V}_t}$, hence $\pardeg \pa{V}^{''} \leq -d_{k+1} = \pardeg \pa{V}'$, and $\rk \pa{V}^{''} < \rk \pa{V}'$, so it is hard not to see that $\pmu(V'') < \pmu(V')$. The condition on the ranks follows immediately from another application of Lemma \ref{lemma:quantitative-two-step}.
\end{proof}

\begin{corollary}\label{cor:stability}
	Let $(C,D)$ be a hyperbolic curve of genus $g$, let $(\pa{V},\nabla)$ be an irreducible flat parabolic
	vector bundle of rank $r < 2\sqrt{g}$ with regular singularities along $D$. Then an analytically general isomonodromic deformation of $\pa{V}$ is stable, and each of its nonzero proper parabolic subbundles $\pa{V'}$ satisfies
	\[\pardeg \pa{V'} \leq -\frac{g - \lfloor r^{2}/4 \rfloor}{r}.\]
\end{corollary}
\begin{proof}
	The parabolic degree of a nonzero proper parabolic subbundle of a semistable isomonodromic deformation is $-d_{k}$ for some $k$, as the quotient by it has nonnegative parabolic degree. Theorem~\ref{thm:stability} rules out every $k$ with $d_{k} < (g - \lfloor r^{2}/4 \rfloor)/r$, and this includes $d_{k} = 0$ since $\lfloor r^{2}/4 \rfloor < g$.
\end{proof}

\printbibliography
\end{document}